\documentclass[tbtags,reqno]{amsart}        

\usepackage[tbtags]{amsmath}

\usepackage{anysize,amsthm,amsopn,amsfonts,amssymb,epsfig}            
\marginsize{2.5cm}{2.5cm}{2.5cm}{2.5cm}

\numberwithin{equation}{section}

\edef\savecatcodeat{\the\catcode`@}
\catcode`\@=11

\def\tb@ifSpecChars#1#2{#1}
\def\tb@ifNoSpecChars#1#2{#2}

\def\tableau{%
  \bgroup% matched in \tb@tableauD
  \@ifstar{\let\Tif\tb@ifNoSpecChars\tb@tableauB}% *, don't use special chars
          {\let\Tif\tb@ifSpecChars\tb@tableauB}}% no *, use special chars

\def\tb@tableauB{% add [] if no [options]
  \@ifnextchar[{\tb@tableauC}{\tb@tableauC[]}}

\def\tb@tableauC[#1]{\hbox\bgroup%
    \let\\=\cr% end line
    \def\bl{\global\let\tbcellF\tb@cellNF}%
    \def\tf{\global\let\tbcellF\tb@cellH}% highlighted cell
    \dimen2=\ht\strutbox \advance\dimen2 by\dp\strutbox%
    \ifx\baselinestretch\undefined\relax%
    \else%
       \dimen0=100sp \dimen0=\baselinestretch\dimen0%
       \dimen2=100\dimen2 \divide\dimen2 by\dimen0%
    \fi%
    \let\tpos\tb@vcenter% default position
    \tb@initYoung% default tableau type
    \tb@options#1\eoo% parse options
    \let\arrow\tb@arrow%
    \dimen0=\Tscale\dimen2%
    \dimen1=\dimen0 \advance\dimen1 by \tb@fframe%
    \lineskip=0pt\baselineskip=0pt% line spacing will be from \vbox to \dimen0
    \def\tb@nothing{}%
    \def\endcellno{$\rss\egroup\bss\egroup}% end cell w/o overlap
    \def\endcell{\endcellno\kern-\dimen0}% end cell & prepare to overlap it
    \def\begincell{\vbox to\dimen0\bgroup\vss\hbox to\dimen0\bgroup\hss$}%
    \let\overlay\tb@overlay%
    \let\fl\tb@fl%
    \let\lss\hss\let\rss\hss\let\tss\vss\let\bss\vss% cell alignment
    \def\mkcell##1{% format individual cell
        \let\tbcellF\tb@cellD% default cell frame
        \def\tb@cellarg{##1}% store cell contents
        \ifx\tb@cellarg\tb@nothing\let\tb@cellarg\tb@cellE\fi%
        \begincell\tb@cellarg\endcellno% the actual cell content
        \tbcellF}% draw cell frame
    \let\savecellF\tbcellF% save global value of cellF in case of nested tableau
     \Tif{\catcode`,=4\catcode`|=\active}{}\tb@tableauD}%

\let\tb@savetableauD\tableauD% save any current definition
{% set up characters which will be interpreted as command characters
    \catcode`|=\active \catcode`*=\active \catcode`~=\active%
    \catcode`@=\active% command characters
\gdef\tableauD#1{%
  \Tif{% make all the command characters active in math mode when #1 parsed
    \mathcode`|="8000 \mathcode`*="8000%
    \mathcode`~="8000 \mathcode`@="8000%
    \def@{\bullet}%
    \let|\cr% end line
    \let*\tf% highlighted cell
    \let~\sk% skew cell
  }{}%
  \tpos{\tabskip=0pt\halign{&\mkcell{##}\cr#1\crcr}}%
  \global\let\tbcellF\savecellF% restore global value
  \egroup% match \hbox\bgroup at start of \tableauC
  \egroup}% match \bgroup at start of \tableau
}
\let\tb@tableauD\tableauD% rename the command
\let\tableauD\tb@savetableauD% restore old command with this name
\let\tb@savetableauD\undefined

\def\tb@options#1{\ifx#1\eoo\relax\else\tb@option#1\expandafter\tb@options\fi}

\def\tb@option#1{%
  \if#1t\let\tpos\tb@vtop\fi%        t = align at top
  \if#1c\let\tpos\tb@vcenter\fi%     c = align at center
  \if#1b\let\tpos\vbox\fi%           b = align at bottom
  \if#1F\tb@initFerrers\fi%          F = Ferrers diagram
  \if#1Y\tb@initYoung\fi%            Y = Young diagram
  \if#1s\tb@initSmall\fi%            s = small boxes
  \if#1m\tb@initMedium\fi%           m = medium boxes
  \if#1l\tb@initLarge\fi%            l = large boxes
  \if#1p\tb@initPartition\fi%            p = small partition sized boxes
  \if#1a\tb@initArrow\fi%            a = use arrow font as base dimension
}

\def\tb@vcenter#1{\ifmmode\vcenter{#1}\else$\vcenter{#1}$\fi}

\def\tb@vtop#1{\hbox{\raise\ht\strutbox\hbox{\lower\dimen0\vtop{#1}}}}

\def\tb@initPartition{\def\Tscale{.3}}
\def\tb@initSmall{\def\Tscale{1}}
\def\tb@initMedium{\def\Tscale{2}}
\def\tb@initLarge{\def\Tscale{3}}

\def\tb@initArrow{\dimen2=1.25em}

\def\tb@initYoung{%
  \def\tb@cellE{}% empty cell stays empty
  \let\tb@cellD\tb@cellN% default frame is normal frame
  \def\sk{\global\let\tbcellF\tb@cellNF}}% skew cells are empty
\def\tb@initFerrers{%
  \def\tb@cellE{\bullet}% empty cell gets bullet
  \let\tb@cellD\tb@cellNF% default frame is no frame
  \def\sk{\bullet}}% skew cell gets bullet

\tb@initMedium% default scale

\def\tb@sframe#1{%
  \vbox to0pt{%            Embed frame in a box of no vert or hor extent
    \vss%                            pull box above reference point
    \hbox to0pt{%
      \hss%                          pull box left of reference point
      \vbox to\dimen1{%              Actual width of frame
        \hrule depth #1 height0pt% draw top edge of frame
        \vss%                     begin vcenter sides
        \hbox to\dimen1{%           horiz box with side edges just inside
          \vrule width #1 height\dimen1% left edge
          \hss%                     stretch center
          \vrule width #1%         right edge
          }%
        \vss%                     end vcenter sides
        \hrule height #1 depth 0in% bottom edge
        }%
      \kern-\tb@hframe%           horiz alignment off by half line width
      }%
    \kern-\tb@hframe}}%           vert alignment off by half line width
\def\tb@hframe{.2pt}\def\tb@fframe{.4pt}\def\tb@bframe{2pt}
\def\tb@cellH{\tb@sframe{\tb@bframe}}       % bold frame
\def\tb@cellNF{}                            % no frame
\def\tb@cellN{\tb@sframe{\tb@fframe}}       % normal frame
\let\tbcellF\tb@cellN                       % default is normal

\def\tb@rpad{1pt}
\def\tb@lpad{1pt}
\def\tb@tpad{1.8pt}
\def\tb@bpad{1.8pt}

\def\tb@overlay{\endcell\@ifnextchar[{\tb@overlaya}{\begincell}}
\def\tb@overlaya[#1]{\vbox to\dimen0\bgroup%
  \tb@overlayoptions#1\eoo%
  \tss\hbox to\dimen0\bgroup\lss$}
\def\tb@overlayoptions#1{\ifx#1\eoo\relax\else\tb@overlayoption#1\expandafter\tb@overlayoptions\fi}

\def\tb@overlayoption#1{
  \if#1t\def\tss{\vskip\tb@tpad}\let\bss\vss\fi% t = align at top
  \if#1c\let\tss\vss\let\bss\vss\fi%             c = align at center
  \if#1b\def\bss{\vskip\tb@bpad}\let\tss\vss\fi% b = align at bottom
  \if#1l\def\lss{\hskip\tb@lpad}\let\rss\hss\fi% l = align at left
  \if#1m\let\lss\hss\let\rss\hss\fi%             m = align at middle
  \if#1r\def\rss{\hskip\tb@rpad}\let\lss\hss\fi% r = align at right
}

\def\tb@fl{\endcell\begincell\vrule depth 0pt width \dimen0 height \dimen0 \endcell\begincell}

\@ifundefined{diagram}{}{
\def\tb@arrowpad{.5}

\newoptcommand{\tb@arrow}{\@ne}[2]{%
  \endcell% end previous cell contents
   \begingroup%
   \let\dg@getnodesize\tb@getnodesize% substitute routine to get nodesize
   \dg@USERSIZE=#1\relax%
   \ifnum\dg@USERSIZE<\@ne \dg@USERSIZE=\@ne \fi%
   \dg@parse{#2}%
   \dg@label{\tb@draw{#1}{#2}}}% draw arrow

\def\tb@getnodesize#1#2#3#4#5{\dimen3=\tb@arrowpad\dimen2 #4=\dimen3 #5=\dimen3\relax}
\def\tb@getnodesize#1#2#3#4#5{\ifnum#2=0\ifnum#3=0\tb@getnodesizetail{#4}{#5}\else\tb@getnodesizehead{#4}{#5}\fi\else\tb@getnodesizehead{#4}{#5}\fi}
\def\tb@getnodesizetail#1#2{\dimen3=.5\dimen2 #1=\dimen3 #2=\dimen3}
\def\tb@getnodesizehead#1#2{\dimen3=.5\dimen2 #1=\dimen3 #2=\dimen3}

\def\tb@draw#1#2#3#4{%
        \dg@X=0\dg@Y=0\dg@XGRID=1\dg@YGRID=1\unitlength=.001\dimen0%
        \dg@LBLOFF=\dgLABELOFFSET \divide\dg@LBLOFF\unitlength%
        \dg@drawcalc% compute arrow geometry
        \begincell% start tableau cell
        \let\lams@arrow\tb@lams@arrow% substitute routine
        \begin{picture}(0,0)\begingroup\dg@draw{#1}{#2}{#3}{#4}\end{picture}%
        \endcell% end tableau cell
        \endgroup% match \begingroup in \tb@arrow
        \begincell}% start new entry in this cell
}

\def\tb@lams@arrow#1#2{%
 \lams@firstx\z@\lams@firsty\z@
 \lams@lastx#1\relax\lams@lasty#2\relax
 \lams@center\z@
 \N@false\E@false\H@false\V@false
 \ifdim\lams@lastx>\z@\E@true\fi
 \ifdim\lams@lastx=\z@\V@true\fi
 \ifdim\lams@lasty>\z@\N@true\fi
 \ifdim\lams@lasty=\z@\H@true\fi
 \NESW@false
 \ifN@\ifE@\NESW@true\fi\else\ifE@\else\NESW@true\fi\fi
 \ifH@\else\ifV@\else
  \lams@slope
  \ifnum\lams@tani>\lams@tanii
   \lams@ht\ten@\p@\lams@wd\ten@\p@
   \multiply\lams@wd\lams@tanii\divide\lams@wd\lams@tani
  \else
   \lams@wd\ten@\p@\lams@ht\ten@\p@
   \divide\lams@ht\lams@tanii\multiply\lams@ht\lams@tani
  \fi
 \fi\fi
 \ifH@  \lams@harrow
 \else\ifV@ \lams@varrow
 \else \lams@darrow
 \fi\fi
}

\catcode`\@=\savecatcodeat
\let\savecatcodeat\undefined

\newtheorem{theorem}{Theorem}

\newtheorem{proposition}[theorem]{Proposition}
\newtheorem{example}[theorem]{Example}
\newtheorem{conjecture}[theorem]{Conjecture}

\newtheorem{corollary}[theorem]{Corollary}
\newtheorem{definition}[theorem]{Definition}
\newtheorem{property}[theorem]{Property}
\newtheorem{remark}[theorem]{Remark}

\def\cpreceq{\prec  \!  \!\!\cdot \,}
\newcommand{\lm}{\lambda/\mu}

\def\H{{\mathcal H}}
\def\sap{\bigskip}
\def\sas{\smallskip}
\def\sa{\medskip}
\def\charge{ {\rm {charge}}}
\def\cocharge{ {\rm {cocharge}}}
\def\shape{ {\rm {shape}}}
\def\stat { {\rm {Stat}}}
\def\dv { {\rm {dv}}}
\def\tab  { \mathsf T }
\def\sstab  {  T }
\def\stab  { \mathcal T }
\def \Sym {\mathcal Sym}
\def \SS {\lambda/^{k}}

\begin{document}

\thispagestyle{plain}

\title[Order ideals in  weak subposets of Young's lattice]
{Order ideals in weak subposets of Young's lattice
and associated unimodality conjectures}

\author{L. Lapointe}
\thanks{Research supported in part by the
Fondo Nacional de Desarrollo Cient\'{\i}fico y Tecnol\'ogico
(FONDECYT) project \#1030114, the Programa Formas Cuadr\'aticas of
the Universidad de Talca, and by NSERC grant \#250904}
\address{
Instituto de Matem\'atica y F\'{\i}sica,
Universidad de Talca,
Casilla 747, Talca, Chile}
\email{lapointe@inst-mat.utalca.cl}

\author{J. Morse}
\thanks{Research supported in part by NSF grant \#DMS-0400628}
\address
{Department of Mathematics,
University of Miami, 
Coral Gables, FL 33124} 
\email {morsej@math.miami.edu}

\begin{abstract}
The $k$-Young lattice $Y^k$ is a weak subposet of the Young lattice
containing partitions whose first part is bounded by an integer $k>0$.
The $Y^k$ poset was introduced in connection with generalized Schur 
functions and later shown to be isomorphic to the weak order on the 
quotient of the affine symmetric group $\tilde S_{k+1}$ by a maximal 
parabolic subgroup.  We prove a number of properties for $Y^k$ including 
that the covering 
relation is preserved when elements are translated by rectangular 
partitions with hook-length $k$.  We highlight the order ideal 
generated by an $m\times n$ rectangular shape.  This order ideal, 
$L^k(m,n)$, reduces to $L(m,n)$ for large $k$, and we prove it is
isomorphic to the induced subposet of $L(m,n)$ whose vertex set 
is restricted to elements with no more than $k-m+1$ parts 
smaller than $m$.  We provide explicit formulas for the number of 
elements and the rank-generating function of $L^k(m,n)$. We 
conclude with unimodality conjectures involving $q$-binomial 
coefficients and discuss how implications connect to recent work
on sieved $q$-binomial coefficients.
\end{abstract}


\maketitle

%%%% Authors begin text of article here %%%

\def\H{{\mathcal H}} 
\def\sap{\bigskip}
\def\sas{\smallskip}
\def\sa{\medskip}
\def\charge{ {\rm {charge}}}
\def\cocharge{ {\rm {cocharge}}}
\def\shape{ {\rm {shape}}}
\def\stat { {\rm {Stat}}}
\def\dv { {\rm {dv}}}
\def\tab  { \mathsf T }
\def\sstab  {  T }
\def\stab  { \mathcal T }
\def \Sym {\mathcal Sym}
\def\endprf {\square}
\def\con {\overline{\in}}
 
\section{Introduction}
The Young lattice $Y$ is the poset of integer partitions given by inclusion 
of diagrams.  This poset can be induced from the branching rules of the 
symmetric group, and certain order ideals of $Y$ are in themselves
interesting posets.  For example, the induced subposet of partitions 
whose Ferrers diagrams fit inside an $m\times n$ rectangle satisfies 
many beautiful properties.  These principal order ideals, denoted
$L(m,n)$, are graded, self-dual, and strongly sperner lattices \cite{[St]}.  
Further, it is known that the number of elements $p_i(m,n)$
of rank $i$ in $L(m,n)$ 
are coefficients in the generalized 
Gaussian polynomial, and thus form a unimodal sequence \cite{[P],[S]}.  
That is,
\begin{equation}
\sum_{i\geq 0}p_i(m,n)\,q^i = 
\begin{bmatrix}
n+m  \\ m
\end{bmatrix}_q = \frac{(1-q^{n+1})\cdots(1-q^{m+n}) }
{(1-q)\cdots(1-q^{m})}
\, .
\label{qbin}
\end{equation}
Letting $q\to 1$, 
the total number of elements in this poset is given by
\begin{equation}
\left|L(m,n) \right|
 = {n+m\choose m}
\, .
\label{bin}
\end{equation}

A weak subposet $Y^k$ of the Young lattice was introduced in 
connection with functions that generalize the Schur functions
\cite{[LLM],[LM1]}.  This poset (hereafter called the $k$-Young lattice) 
is a lattice defined on the set of partitions whose first part is no larger than 
fixed $k\geq 1$.  The order arises from a degree preserving 
involution on the set of $k$-bounded partitions $\mathcal P^k$
that generalizes partition conjugation.  The involution 
sends one $k$-bounded partition $\lambda$ to another, 
$\lambda^{\omega_k}$, giving rise to a partial order on 
$\mathcal P^k$ as follows:
For $\lambda$ and $\mu$ differing by one box,
\begin{quote}
\it
Young order on partitions:
$\lambda<\!\!\cdot\,\mu$ when $\lambda\subseteq\mu$ 
(and equivalently $\lambda'\subseteq\mu'$).
\end{quote}
\begin{quote}
\it
$k$-order on $k$-bounded partitions: $\lambda\prec\!\!\cdot\,\mu$ when 
$\lambda\subseteq\mu$ and $\lambda^{\omega_k}\subseteq\mu^{\omega_k}$.
\end{quote}
It happens that $\lambda^{\omega_k}=\lambda'$
for large $k$ implying that the $k$-order is the Young order
in the limit  $k\to \infty$.

The $k$-Young lattice originated from a conjectured formula for 
multiplying $k$-Schur functions \cite{[LLM],[LM1]} that is analogous to the 
Pieri rule.  In particular, the conjecture states that the $k$-Schur functions 
$s_{\mu}^{(k)}$ appearing in the expansion of the
product $s_1 s_{\lambda}^{(k)}$ are exactly those indexed by 
the successors of $\lambda$ in the $k$-Young lattice.  That is,
$$
s_1 s_{\lambda}^{(k)}=\sum_{\lambda\prec\!\!\cdot\mu} s_{\mu}^{(k)}\, .
$$
In \cite{[LM3]}, it is shown that the $k$-Young lattice is in fact
isomorphic to the weak order on the quotient of the affine symmetric 
group $\tilde S_{k+1}$ by a maximal parabolic subgroup and that
the paths in $Y^k$ can be enumerated by certain ``$k$-tableaux", 
or by reduced words for affine permutations.

Here we investigate general properties of the $k$-Young lattice.  Most notably,
we reveal that partitions with rectangular shape and hook-length $k$, 
called $k$-rectangles, play a fundamental role in the structure of this 
poset.  We prove for any $k$-rectangle $\square$ and
$\lambda,\mu\in\mathcal P^k$,
$$ \lambda \preceq \mu\;\;\text{  if and only if }\;\;
(\lambda\cup\square) \preceq (\mu \cup \square)\,,
$$
leading to the stronger statement that:
$$
\lambda\, \cup \,\square \preceq\mu\iff
\mu = \bar \mu \, \cup \, \square \quad and \quad \lambda \preceq \bar \mu\,,
$$
for some $k$-bounded partition $\bar \mu$.  This is a central property 
needed to identify the $k$-Young lattice with a cone in the 
permutahedron-tiling of the $k$-dimensional space \cite{[Ul]}. 
The significance of $k$-rectangles also plays an important role at the
symmetric function level in that multiplying a Schur function indexed 
by a $k$-rectangle $\square$ with a $k$-Schur function is trivial 
\cite{[LMrec]}.  That is, for any $k$-bounded partition $\lambda$, 
$$
s_{\square} \, s_{\lambda}^{(k)}=s_{\square \, \cup \, \lambda}^{(k)}
\,.
$$

Following our study of the $k$-rectangles and other properties of
the $k$-Young lattice, we discuss a
family of induced subposets of $L(m,n)$ whose vertex set 
consists of the elements that fit inside an $m\times n$ rectangle 
and have no more than $k-m+1$ parts strictly smaller than $m$.
Surprisingly, we find that these subposets are isomorphic to the 
principal order ideal of $Y^k$ generated by the shape $m\times n$.
As such, we denote these order ideals by $L^k(m,n)$ and
note that they are graded, self-dual, distributive lattices
of rank $mn$.
We provide explicit formulas for the number of vertices and the 
rank generating function:  
\begin{equation}
\sum_{\lambda\in L^{k}(m,n)} 
q^{|\lambda|}
\,=\,
\begin{bmatrix}
k+1 \\m
\end{bmatrix}_q
+
q^{k+1}
\frac{1-q^{m(n-k+m-1)}}{1-q^m}
\begin{bmatrix}
k \\m-1
\end{bmatrix}_q
\, ,
\label{eqbink}
\end{equation}
for $n \geq k-m+1$,
which implies that the coefficient of $q^i$ in the right hand 
side of this expression is the number of partitions in $L^k(m,n)$ 
with rank $i$.  Further, letting $q\to 1$, 
\begin{equation}
\left|L^k(m,n) \right|
 = 
\binom{k+1}{m}+(n-k+m-1)
\binom{k}{m-1} 
\, ,
\label{krec1}
\end{equation}
for $n\geq k-m+1$.

Since the vertex set of $L^k(m,n)$ is contained in that of
$L^{k+1}(m,n)$, these order ideals provide a natural sequence of
subposets of $L(m,n)$.  That is, $L(m,n)$ can be constructed from
the chain of partitions with no more than one row smaller 
than $m$ by successively adding sets of partitions with
exactly $j$ parts smaller than $m$ for $j\geq 2$.
This decomposition aids our investigation of questions
pertaining to unimodality.  Prompted by the unimodality 
of $L(m,n)$, we computed examples that suggest $L^k(m,n)$ is 
unimodal in certain cases. In particular when $k\neq -1\mod p$ 
for all prime divisor $p$ of $m$.  When $m$ is prime, we find that
the unimodality of $L^k(m,n)$ relies on the conjecture:
\begin{quote}
{\it
If $k\neq -1,0 \mod m$,
then the coefficients of the $q$ powers in
\begin{equation}
\label{introuni}
\frac{(1-q^{m(n-k+m)})}{(1-q^m)} 
\begin{bmatrix}
k-1 \\m-2
\end{bmatrix}_q
\end{equation}
form a unimodal sequence for all $n \geq k-m+1$. 
}
\end{quote}
We also generalize this conjecture to include the case
when $m$ is not prime (see Conjecture~\ref{conjecu}).  

We conclude with a discussion of how our conjectures lead to 
results coinciding with recent work on sieved binomial 
polynomials, eg. \cite{[GS],[SW],[WW]}.  Namely, from the 
unimodality of the coefficients in Eq.~\eqref{introuni}, we 
recover the identity:
the sum of the coefficients of $q^{\ell+*m}$ 
in $ 
\begin{bmatrix}
k-1 \\m-2
\end{bmatrix}_q$
is equal to
$\frac{1}{m}\begin{pmatrix}
k-1 \\m-2
\end{pmatrix}$
if $k \neq -1,0 \mod m$ 
for a prime $m<k$. 
Similarly, we use our more general conjecture to suggest 
a new identity of this type and provide an independent proof 
(see Proposition~\ref{conj4}).

\section{Definitions}

For definitions and general properties of posets see for example 
\cite{[StEnu]}.  A partition $\lambda=(\lambda_1,\dots,\lambda_m)$ is a 
non-increasing sequence of positive
integers. We denote by $|\lambda|$ the degree $\lambda_1 +\cdots +\lambda_m$
of $\lambda$, and by $\ell(\lambda)$ its length $m$.
 Each partition $\lambda$ has an associated Ferrers diagram
with $\lambda_i$ lattice squares in the $i^{th}$ row,
from the bottom to top.  For example,
\begin{equation}
\lambda\,=\,(4,2)\,=\,
{\tiny{\tableau*[scY]{ & \cr & & & }}} \, .
\end{equation}
A partition $\lambda$ is $k$-{\bf bounded} if
$\lambda_1 \leq k$ and the set of all $k$-bounded partitions
is denoted $\mathcal P^k$.  For partitions $\lambda$ and $\mu$,
the weakly decreasing rearrangement of their parts is denoted
$\lambda \cup \mu$, while $\lambda+\mu$ is the partition obtained 
by summing their respective parts.  Any lattice square in the Ferrers 
diagram is called a cell, where the cell $(i,j)$ is in the $i$th row 
and $j$th column in the diagram. We say that $\mu\subseteq \lambda$ when 
$\mu_i\leq \lambda_i$ for all $i$.  

When $\mu \subseteq \lambda$, the skew shape $\lm$ is identified with its
diagram $\{(i,j) : \mu_i<j\leq \lambda_i\}$.  For example,
\begin{equation}
\label{dumbskew}
\lm\,=\,(5,5,4,1)/(4,2)\,=\,
{\tiny{\tableau*[scY]
{\cr&&\tf&\cr\bl&\bl&\tf&&\cr\bl &\bl &\bl s &\bl & \tf}}} \, .
\end{equation}
Lattice squares that do not lie inside a diagram will simply
be called squares.  We shall say that any $s\in\mu$ 
lies {\bf below} the diagram of $\lm$.  The degree of a skew-shape 
is the number of cells in its diagram.
Associated to $\lm$, the {\bf hook} of any $(i,j)\in\lambda$
is defined by the cells of $\lm$ that lie in the L formed with
$(i,j)$ as its corner.  This definition is well-defined for
all squares in $\lambda$ including those below $\lm$.
For example, the framed cells in \eqref{dumbskew} denote the hook of square 
$s=(1,3)$.
We then let $h_s(\lm)$ denote the hook-length of any
$s\in\lambda$, i.e. the number of cells in the hook of $s$.
For example,
$h_{(1,3)}(5,5,4,1)/(4,2)=3$ and $h_{(3,2)}(5,5,4,1)/(4,2)=3$
or cell (3,2) has a 3-hook.  We also say that the hook of a cell
(or a square) is $k$-bounded if it is not larger than $k$.

A {\bf removable corner} is a cell $(i,j)\in\lm$ 
such that $(i+1,j), (i,j+1)\not\in \lm$,
and an {\bf addable corner}  is a square $(i,j)\not\in\lm$ such
that  $(i-1,j), (i,j-1)\in\lm$.  We shall include
$(1,\lambda_1)\in\lambda/\mu$ as a removable corner, 
and $(\ell(\lambda)+1,1)$ as an addable corner.
The {\bf $k$-residue} of any cell (or square) $(i,j)$ in a skew-shape $\lm$
is $j-i \mod k$. That is, the integer in this cell (or square) when $\lm$ 
is periodically labeled with $0,1,\ldots,k-1$, where zeros fill the 
main diagonal.  For example, cell $(1,5)$ has 5-residue 4:
$$
{\tiny{\tableau[scY]{\bl 1|2,\bl 3|\bl 3,4,0,\bl 1|\bl 4,0,1,2,
\bl 3|\bl 0,\bl 1,\bl 2,3,4,0,\bl 1}}}
\, .
$$

\begin{remark}
\label{recres}
No two removable corners of any partition fitting inside the 
shape $\left(m^{k-m+1}\right)$ have the same $k+1$-residue
for any $1\leq m\leq k$.
\end{remark}

\section{Involution on $k$-bounded partitions}

Usual partition conjugation defined by the column reading of
diagrams does not send the set of $k$-bounded partitions to itself.  
Thus our study of $\mathcal P^k$ begins with the 
need for a degree-preserving involution on this set that
extends the notion of conjugation.  Such an involution on 
$\mathcal P^k$ was defined in \cite{[LLM]} using a certain 
subset of skew-diagrams.  We shall follow the notation of 
\cite{[LM3]}, where these skew diagrams are defined by:

\begin{definition} 
The {\bf $k$-skew diagram} of a $k$-bounded partition $\lambda$
is the skew diagram, denoted $\lambda/^k$, satisfying the conditions:

\noindent
(i) row $i$ of $\lambda/^k$ has length $\lambda_i$

\noindent
(ii) no cell in $\lambda/^k$ has a hook-length exceeding $k$

\noindent
(iii) every square below the diagram of $\lambda/^k$
has hook-length exceeding $k$
\label{kskew}
\end{definition}

\begin{example}
Given $\lambda =(4,3,2,2,1,1)$ and $k=4$,
\begin{equation*}
\lambda = {\tiny{\tableau*[scY]{  \cr  \cr  & \cr &
\cr & & \cr & & & }}}
\quad
\implies
\quad
\lambda/^4 = {\tiny{\tableau*[scY]{  
\cr  
\cr  & 
\cr \bl &  & 
\cr  \bl &\bl & & & 
\cr \bl &\bl& \bl & \bl &\bl & & & &
}}}
\end{equation*}
\label{exskew}
\end{example}

It was shown in \cite{[LM3]} that $\lambda/^k$
is the unique skew diagram obtained recursively by:
\begin{quote} 
For any $\lambda=(\lambda_1,\ldots,\lambda_\ell)\in\mathcal P^k$, 
$\lambda/^k$ can be obtained by adding to the bottom of 
$\left(\lambda_2,\ldots,\lambda_\ell\right)/^k$,
a row of $\lambda_1$ cells whose first ({\it i.e.} leftmost) cell $s$  
occurs in the leftmost column where $h_s\leq k$.
That is, row  $\lambda_1$ lies as far to the left as possible
without violating Condition (ii) of $\lambda/^k$, 
or without creating a non-skew diagram.
\end{quote}
As a matter of curiosity, with the skew diagram $\lambda/^k=\gamma/\rho$, 
it is shown in \cite{[LM3]} that
a bijection between $k$-bounded partitions and $k+1$-cores arises
by taking $\lambda\to \gamma$.

Since the columns of a $k$-skew diagram form a partition \cite{[LM3]},
and the transpose of a $k$-skew diagram clearly satisfies Conditions~(ii) 
and (iii), an involution on $\mathcal P^k$ arises from 
the column reading of $\lambda/^k$:

\begin{definition}
For any $k$-bounded partition $\lambda$,
the {\bf $k$-conjugate} of $\lambda$ is the partition given by 
the columns of $\lambda/^k$ and denoted $\lambda^{\omega_k}$.  Equivalently,
$\lambda^{\omega_k}$ is the unique $k$-bounded partition such that $(\lambda/^k)'=\lambda^{\omega_k}/^k$.
\label{defconju}
\end{definition}

\begin{corollary}
For a $k$-bounded partition $\lambda$, we have $(\lambda^{\omega_k})^{\omega_k}=\lambda$.
\end{corollary}

\begin{example}
With $\lambda$ as in Example~\ref{exskew},
the columns of $\lambda/^4$ give
$\lambda^{\omega_4} = (3,2,2,1,1,1,1,1,1)$.
\end{example}

Given a partition $\lambda$, the nature of its $k$-conjugate
is not revealed explicitly by Definition~\ref{defconju}.  However,
the $k$-conjugate can be given explicitly in certain cases.
For example,

\begin{remark}
If $h_{(1,1)}(\lambda)\leq k$, all hooks of $\lambda$
are $k$-bounded and thus $\lambda/^k=\lambda$.
In this case, $\lambda^{\omega_k}=\lambda'$.
\label{re2}
\end{remark}
We can also give a formula for the $k$-conjugate
of a partition with rectangular shape.

\begin{proposition}
$(m^n)^{\omega_k}=\bigl((k-m+1)^a,b^m \bigr)$ where
$b=n\!\mod k-m+1$
and $a=m \lfloor \frac{n}{k-m+1} \rfloor$.
\label{kconrec}
\end{proposition}

\begin{proof}
Building the $k$-skew diagram recursively from the
partition $(m^n)$ reveals that the top $k-m+1$ rows 
are stacked in the shape of the rectangle $\square=(m^{k-m+1})$. 
However, the $k-m+2$-{nd} row (of size $m$) cannot lie below any 
column of $\square$ without creating a $k+1$-hook since
all the columns of $\square$ have height $k-m+1$.  
Thus, it lies strictly to the right of $\square$.
By iteration, $(m^n)/^k$ is comprised of a sequence of block 
diagonal rectangles $\square$ followed by a rectangular block of 
size $(m^{n\mod k-m+1})$.  For example,
\begin{equation}
(3^7) = {\tiny{\tableau*[scY]{  && \cr  && \cr & & \cr  & & \cr & & \cr  & & \cr & & \cr}}} \quad
\longrightarrow
(3^7)/^4 = 
{\tiny{\tableau*[scY]{  && \cr  && 
\cr \bl & \bl & \bl & & &
\cr \bl & \bl & \bl & & &
\cr \bl & \bl & \bl & \bl & \bl & \bl & & &  
\cr \bl & \bl & \bl & \bl & \bl & \bl & & &  
\cr \bl & \bl & \bl & \bl & \bl & \bl & \bl \bl & \bl & \bl & & &   
}}}
\, .
\end{equation}
The columns of such a skew-shape are thus as indicated.
\end{proof}

The rectangular blocks occurring in the $k$-skew diagram of
$(m^n)$ have the form $(m^{k-m+1})$ for $1\leq m\leq k$.   
We have found that such rectangles, called {\bf $k$-rectangles}, 
play an important role in our study.  For starters, we show
that $k$-conjugation can be distributed over the union of any 
partition and a $k$-rectangle. To prove this result,
we first need to find the $k$-conjugate of another shape:

\begin{proposition}
If $\square=(\ell^{k-\ell+1})$ and 
$\mu\subseteq (\ell-1)^{k-\ell}$, then
\begin{equation}
(\square,\mu)/^k=(\square+\mu,\mu)/\mu
\quad\text{implying}
\quad (\square,\mu)^{\omega_k}=(\square',\mu')\,.
\end{equation}
\label{recskew2}
\end{proposition}
\begin{proof}
Given $\mu \subseteq(\ell-1)^{k-\ell}\subseteq \square$, 
$D=(\square+\mu,\mu)/\mu$ is a skew diagram with 
general shape depicted in Figure~\ref{nouveau}.
\begin{figure}[htb]
\epsfig{file=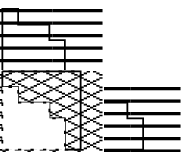}
\caption{$(\square+\mu,\mu)/\mu$, with $\mu$ depicted 
by the partition in
horizontal stripes}
\label{nouveau}
\end{figure}
If $D$ meets Conditions~(ii) and (iii) for a $k$-skew diagram 
then \hbox{$(\square,\mu)/^k=(\square+\mu,\mu)/\mu$} since the rows of 
$D$ are given by the partition $(\square,\mu)$ and thus
$(\square,\mu)^{\omega_k}=(\square',\mu')$
since the columns of $D$ are $(\square',\mu')$.
Any square below $D$ has $k-\ell+1$ cells above it and $\ell$ to the 
right implying it has hook $k+1>k$. On the other hand, any cell in 
$D$ has hook-length strictly smaller than this since the columns 
and rows
of $D$ are weakly decreasing.  Therefore $D=(\square,\mu)/^k$.
\end{proof}

\begin{theorem} \label{rectangle}
$(\lambda \cup \square)^{\omega_k}= 
\lambda^{\omega_k}\cup \square^{\omega_k}$
for any $k$-bounded partition $\lambda$ and $k$-rectangle $\square$.
\end{theorem}
\begin{proof}
Let $i$ be such that $\lambda_i<\ell$ and $\lambda_{i-1} \geq \ell$,
and let $\mu$ denote the non-skew partition determined by 
the cells strictly above the bottom row of
$(\ell,\lambda_i,\dots,\lambda_{\ell(\lambda)})/^k$ (Figure~\ref{upper}). 
\begin{figure}[htb]
\begin{center}
\epsfig{file=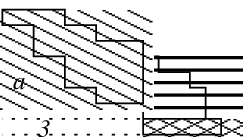}
\end{center}
\caption{$(\ell,\lambda_i,\dots,\lambda_{\ell(\lambda)})/^k$
with $\mu$ depicted by the partition in horizontal stripes}
\label{upper}
\end{figure}
The squares in regions $(a)$ and $(3)$ in the picture
have hooks exceeding $k$ by definition of $k$-skews.  
Now consider Figure~\ref{tout}
\begin{figure}[htb]
\begin{center}
\epsfig{file=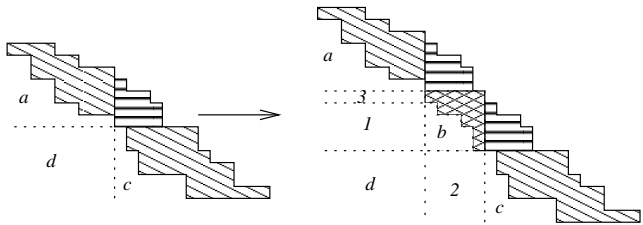}
\end{center}
\caption{Comparison of $\lambda/^k$ and $D$}
\label{tout}
\end{figure}
where the diagram on the left is $\lambda/^k$ and 
in the diagram $D$ on the right, region $(b)$ is of shape $\mu$.
Since the rows of $D$ are $\lambda \cup \square$,
if we prove that $D$ is a $k$-skew diagram, then
$D=(\lambda \cup \square)/^k$.
This will then imply $(\lambda \cup \square)^{\omega_k}= 
\lambda^{\omega_k}\cup \square^{\omega_k}$
since the columns of $D$ are just the columns of 
$\lambda/^k$ and $\square$.

The hooks in regions $(a)$ and $(c)$ of $D$
are the same as they are in $\lambda/^k$, and thus exceed $k$.
Similarly, the hooks to the right and above $(a)$ and $(c)$ 
are $k$-bounded in $D$.
The squares in region $(3)$ lie ``below" the $k$-skew diagram 
$(\ell,\lambda_i,\ldots,)/^k$ and thus exceed $k$,
as do all squares in regions $1$ and $d$ since
they can only increase given that the rows of $D$ form a partition.
The cells of region $(2)$ also exceed $k$ since there are
$k-\ell+1$ cells above them and at least $\ell$ to their right.  

Finally, the subdiagram of $D$ including region $(b)$ and all cells
above and to the right of this region is the diagram of 
$(\square+\mu,\mu)/\mu$.
Since $\mu_1 < \ell$ and $\ell(\mu) \leq k-\ell$ by definition,
Proposition~\ref{recskew2} implies that this subdiagram
is $(\square,\mu)/^k$ and thus meets the conditions of a $k$-skew.
Therefore, $D$ is a $k$-skew diagram and the theorem follows.
\end{proof}

Ideas to understand the nature of a $k$-skew diagram containing
a $k$-rectangle that were used in Theorem~\ref{rectangle} may be applied to 
prove the following technical proposition to be used later.

\begin{proposition} \label{hookegalk}
For some $1\leq\ell\leq k$, if $\nu$ is a partition containing 
exactly $k-\ell+1$ rows of length $\ell$, where the 
lowest occurs in some row $r$, then there are addable corners in row 
$r$ and $k-\ell+1+r$ of  $\nu/^k$ with the same  $k+1$-residue.
\end{proposition}
\begin{proof}
Let $\square=(\ell^{k-\ell+1})$ and $\nu=\lambda \cup \square$ 
for some partition $\lambda$ with no parts of size $\ell$. 
We can construct the diagram of $\nu/^k$ as in the previous proof.
That is, let $i$ be such that $\lambda_i<\ell$ and 
$\lambda_{i-1} > \ell$, and let $\mu$ denote the 
non-skew partition determined by the cells strictly above 
the bottom row of $(\ell,\lambda_i,\dots,\lambda_{\ell(\lambda)})/^k$.  
We appeal to Figure~\ref{tout}, where the 
diagram of $(\lambda \cup \square)/^k$ is on the right.

Note that if $r$ denotes the lowest row in $\nu/^k$ of length $\ell$,
then row $r-1$ (if it exists) is strictly longer than
row $r$ since $\lambda_{i-1}>\ell$.
Therefore an addable corner $x$ occurs in row $r$ with some
$k+1$-residue $j$.  Furthermore, since row $k-\ell+1+r$ 
corresponds to row $\lambda_{i-1}$ (the first row of $\mu$) 
and $\lambda_{i-1}<\ell$,
there is also an addable corner $\bar x$ in this row.
Thus, it remains to show that $\bar x$ has $k+1$-residue $j$.  
Let $s$ denote the first cell in row $r$.  If $\bar x$ lies 
in the column of $s$, then the hook length of cell $s$ is  
$k-\ell+\ell=k$, and thus $\bar x$ also has $k+1$-residue $j$.  
We shall now see that $\bar x$ does in fact lie in the column 
containing $s$.  If $\bar x$ lies in a
column to the right of $s$, then since $\bar x$ is an addable corner,
the hook-length of $s$ is
larger than $k$ (a contradiction).  And
if $\bar x$ lies in a column to the left of $s$, then 
the square in the row of $x$ and the column of $\bar x$ 
has hook-length equal to at most $k-\ell+\ell=k \not >k$.
\end{proof}

\section{$k$-Young lattices}

Recall that the Young Lattice $Y$ is the poset of all partitions 
ordered by inclusion of diagrams, or equivalently 
$\lambda\leq \mu$ when $\lambda \subseteq \mu$.
Since $\lambda\subseteq\mu\iff\lambda'\subseteq\mu'$,
it is equivalent to view the Young order as
$\lambda\leq \mu$ when $\lambda \subseteq \mu$ {\it and}
$\lambda'\subseteq\mu'$.  
This interpretation for the Young order refines naturally
to an order on $k$-bounded partitions by using the the 
$k$-conjugate on $\mathcal P^k$.

\begin{definition}  
\label{korder}
The order $\preceq$ on partitions in $\mathcal P^k$
is defined by the transitive closure of the relation
\begin{equation}
\mu\cpreceq\lambda\quad\text{when}
\quad
\lambda \subseteq \mu \quad\text{ and}\quad
\lambda^{\omega_k}\subseteq \mu^{\omega_k}
\;\;\text{and}\;\; |\mu|-|\lambda|=1\,.
\end{equation}
\end{definition}
\begin{figure}[htb]
\epsfig{file=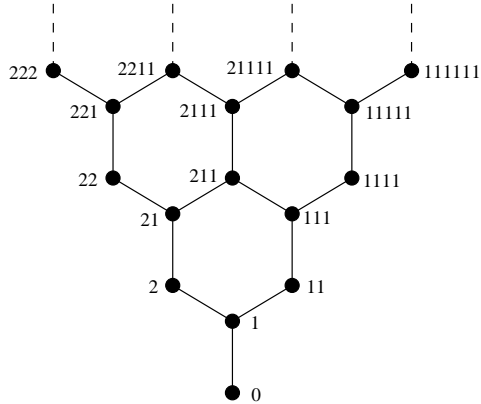}
\caption{Hasse diagram of the $k$-Young lattice in the case $k=2$.}
\label{posetexample}
\end{figure}
We denote this poset on $\mathcal P^k$ by $Y^k$, and observe that
it is a weak subposet of the Young lattice (recall
this means that if $\lambda \preceq \mu$ in $Y^k$, then
$\lambda \preceq \mu$ in $Y$).  Furthermore, $Y^k$ 
reduces to the Young lattice when $k\to \infty$ since
\cite{[LM3]}:

\begin{property}
\label{kord}
$\lambda\preceq\mu$ reduces to $\lambda\leq \mu$
when $\lambda$ and $\mu$ are partitions 
with $h_{(1,1)}(\lambda)\leq k$ and $h_{(1,1)}(\mu)\leq k$.
\end{property}

While this poset on $k$-bounded partitions originally arose in 
connection to a rule for multiplying generalized Schur functions 
\cite{[LLM]}, it has been shown in \cite{[LM3]} that this poset 
turns out to be isomorphic to the weak order on the quotient 
of the affine symmetric group by a maximal parabolic subgroup.  
Consequently, $Y^k$ is a lattice \cite{[Ul],[DW]} and we thus call it the 
{\it $k$-Young lattice}.

Although the ordering $\preceq$ is defined by the covering relation 
$\cpreceq$, it follows from the definition that

\begin{property}
\label{contains}
If $\lambda \preceq \mu$, then $\lambda \subseteq \mu$ and 
$\lambda^{\omega_k}\subseteq \mu^{\omega_k}$.  
\end{property}

It is important to note that the converse of this statement does not hold.
For example, with $k=3$, $\lambda=(2,2)$, and $\mu=(3,2,1,1,1,1)$:
$\lambda^{\omega_k}=\lambda$ and $\mu^{\omega_k}=\mu$
satisfy $\lambda \subseteq \mu$ and
$\lambda^{\omega_k}\subseteq \mu^{\omega_k}$,
but $\lambda \not \preceq \mu$ (see Theorem~\ref{domirect}
and note that $\lambda$ contains the 3-rectangle $(2,2)$ while
$\mu$ does not).

Since the set of $\mu$ such that $\mu\subseteq\lambda$ and 
$|\mu|=|\lambda|-1$ consists of all partitions obtained by 
removing a corner box from $\lambda$, the set of elements 
covered by $\lambda$ with respect to $\preceq$ is a subset 
of these partitions.  The corners that can be removed 
from $\lambda$ to give partitions covered by $\lambda$
are determined as follows:

\begin{theorem}
\cite{[LM3]}
The order $\preceq$ can be characterized by the covering relation
\begin{equation}
\lambda\cpreceq\mu
\quad\iff\quad \mu=\lambda+e_r\, ,
\end{equation}
where $r$ is any row of $\mu/^k$ with a removable corner
whose $k+1$-residue does not occur in a higher removable corner, 
or equivalently for $r$ a row in $\lambda/^k$ with an addable 
corner whose $k+1$-residue does not occur in a higher addable corner.
\label{khook}
\label{khookplus}
\end{theorem}

\begin{example}  \label{whatcorners}
With $k=4$ and $\lambda=(4,2,1,1)$, 
\begin{equation}
\lambda/^4 = {\tiny{\tableau*[scY]{  
\bl 1\cr 
2\cr  
3 & \bl 4\cr  
4 & 0 & \bl 1 \cr 
\bl 0 & \bl 1 & 2 & 3  & 4 & 0 & \bl 1 
}}} \, ,
\end{equation}
and thus the partitions that are covered by $\lambda$ are
$(4,1,1,1)$, and $(4,2,1)$, while those that
cover it are $(4,2,1,1,1)$ and $(4,2,2,1)$.
\end{example}

Since the conditions of Theorem~\ref{khook} are always satisfied 
when choosing the removable corner in the top row of $\mu$ 
we have the corollary:

\begin{corollary}
If $\mu=(\mu_1,\ldots,\mu_\ell)$ is a $k$-bounded
partition, then $\lambda\preceq\mu$ for
$\lambda=(\mu_1,\ldots,\mu_\ell-1)$.
\label{removebox}
\end{corollary}

The conditions of Theorem~\ref{khook} imply that $\mu=\lambda+e_r$.
Therefore:

\begin{property} \label{newcorner} 
\label{newcornerajout} 
Any row of $\mu/^k$ containing a removable corner whose $k+1$-residue 
does not occur in a higher removable corner, corresponds to a row of 
$\mu$ with a removable corner.
\end{property}

As discussed in the introduction, the $k$-rectangles play a fundamental 
role in the study of $k$-Young lattices.  Since $k$-conjugate distributes over 
the union of $k$-rectangles with a partition,
and the $k$-Young lattice relies on $k$-conjugates, we are able to show that
the order is preserved under union with a $k$-rectangle.

\begin{proposition}
\label{recorder} For any $k$-rectangle $\square$,
$\lambda \preceq \mu$  if and only if
$(\lambda\cup\square) \preceq (\mu \cup \square)$.
\end{proposition}
\begin{proof} 
It suffices to consider the case that $\lambda\cpreceq\mu$.  
Given $\lambda \subseteq \mu$ and 
$\lambda^{\omega_k} \subseteq \mu^{\omega_k}$
with $|\mu|-|\lambda|=1$,
clearly $(\lambda\cup \square) \subseteq (\mu \cup \square)$ 
with $|\mu\cup \square|-|\lambda\cup \square|=1$.
Theorem~\ref{rectangle} then implies that
$(\lambda \cup \square)^{\omega_k}=
(\lambda^{\omega_k}\cup \square^{\omega_k})
\subseteq (\mu^{\omega_k} \cup \square^{\omega_k})
= (\mu\cup \square)^{\omega_k}$. That is,
$(\lambda\cup\square) \cpreceq (\mu \cup \square)$.
\end{proof}

In fact, we have a stronger result that amounts to
saying the $k$-rectangles $\square$ play a trivial role when moving 
up in the $k$-Young lattice.  That is, the partitions dominating 
$\lambda \cup \square$ can be obtained by adding the parts of
$\square$ to the partitions that dominate $\lambda$.
The (increasing) covering relations around $\lambda$ and 
$\lambda \cup \square$ are isomorphic, 
and thus these relations are preserved under translation 
by a $k$-rectangle.

\begin{theorem} \label{domirect}
For $\lambda,\mu\in\mathcal P^k$
and a $k$-rectangle $\square$, 
$$
\lambda\, \cup \,\square \preceq\mu\iff
\mu = \bar \mu \, \cup \, \square \quad and \quad \lambda \preceq \bar \mu \, ,$$
for some $k$-bounded partition $\bar \mu$.
\end{theorem}
\begin{proof}
Let the $k$-rectangle be $\square=(\ell^{k-\ell+1})$.
($\Longleftarrow$) follows from Proposition~\ref{recorder}.
For ($\Longrightarrow$), it suffices to consider 
$\lambda\,\cup\,\square\cpreceq\mu$. 
Theorem~\ref{khookplus} implies that $\mu=(\lambda\,\cup\,\square)+e_r$ 
for some row $r$ with an addable corner $o$ whose $k+1$-residue does 
not occur in any higher addable corner of $(\lambda\,\cup\,\square)/^k$.  
Assume by contradiction that 
$(\lambda\,\cup\,\square)+e_r \neq \bar \mu\,\cup\,\square$ 
for any $k$-bounded partition $\bar\mu$.  The only scenario 
where the number of rows of length $\ell$ is reduced by adding 
a box is if $\lambda\,\cup\,\square$ has exactly $k-\ell+1$ rows
of length $\ell$ and row $r$ is the lowest row of length $\ell$.
Thus, by Proposition~\ref{hookegalk}, there is an addable corner
in row $k-\ell+1+r$ of $(\lambda\cup\square)/^{k}$ with the same
$k+1$-residue as $o$.  However, row $k-\ell+1+r$ is higher than
row $r$ and by contradiction,
$\mu= \bar\mu\cup\square$ for some $\bar\mu$.
Finally, given $\lambda \, \cup \, \square \preceq \bar \mu \, \cup 
\, \square$, the 
previous proposition implies $\lambda \preceq \bar \mu$.
\end{proof}

\begin{remark}
Consider $\lambda,\mu\in\mathcal P^k$ and a $k$-rectangle $\square$.
Notice that in general,
$$
\lambda\, \preceq\mu \cup \square\;\; \;
\not\!\!\!\!\!\iff
\lambda = \bar \lambda \, \cup \, 
\square \;\; and \;\; \bar \lambda \preceq \mu \, .
$$ 
For example, with $k=3$: the 3-rectangle 
$(2,1)\preceq(2,2)$ while $(2,1)\not\supseteq(2,2)$.
However,
$$
\square\subseteq\mu
\;\; and \;\;
\lambda\, \preceq\mu \cup \square\;\iff
\lambda = \bar \lambda \, \cup \, 
\square \;\; and \;\; \bar \lambda \preceq \mu \, $$ 
follows from Proposition~\ref{recorder}.  
In this case, what occurs above and below $\mu$ is replicated 
at $\mu \cup \square$.  Interpreting the poset as a cone in a tiling of 
$k$-space by permutahedrons \cite{[Ul]},
this implies that a vertex $\mu$ lying at least a distance
$|\square|$ from the boundary of the cone can not be distinguished
from $\mu\cup\square$, and thus the $k$-rectangles are the vectors of
translation invariance in the tiling.
\end{remark}

\section{Principal order ideal}

Let $Y_\lambda$ denote the principal order ideal generated
by $\lambda$ in the Young lattice.  When $\lambda$ is a 
rectangle, the order ideal is denoted $L(m,n)$ and is the 
induced poset of partitions with at most $n$ parts and largest 
part at most $m$.  This order ideal is a graded, self-dual, and 
distributive lattice.  Further, the rank-generating function of $L(m,n)$, 
the Gausssian polynomial, is known to be unimodal.  The next several 
sections concern the study of properties for order ideals in
the $k$-Young lattice that are analogous to those held by $L(m,n)$.

Let $Y_{\lambda}^k$ denote the principal order ideal generated 
by $\lambda$ in the poset $Y^k$.  That is, 
$$
Y_\lambda^k  = \left\{\mu : \mu\preceq\lambda\right\}\,.
$$
As $k\to\infty$, the poset $Y^k_\lambda$ reduces to $Y_\lambda$.  
More precisely,

\begin{property}
If $\lambda$ is a partition with $h_{(1,1)}(\lambda)\leq k$, then
$Y_\lambda^k=Y_\lambda$.
\label{hookk}
\end{property}

\begin{proof}
Since $\mu\preceq\lambda$ implies that $\mu \subseteq \lambda$, 
we have that $h_{(1,1)}(\mu)\leq 
h_{(1,1)}(\lambda)\leq k$.
Thus, for all $\mu,\nu\in Y_\lambda^k$ we have 
$\mu\preceq\nu\iff\mu\leq\nu$ by
Property~\ref{kord}, that is, $Y_\lambda^k=Y_\lambda$. 
\end{proof}

\begin{proposition}
$Y^k_{\lambda}$ is graded of rank $|\lambda|$.
\label{rank}
\end{proposition}
\begin{proof}
If $\lambda^{(1)}, \lambda^{(2)}, \ldots \lambda^{(n)}$
is a saturated chain in $Y_{\lambda}^k$ then
$|\lambda^{(i)}|= |\lambda^{(i+1)}|-1$ from the definition
of the order $\preceq$.
Corollaries~\ref{removebox} implies that 
a maximal chain in $Y_{\lambda}^k$ must begin with the empty 
partition.  Therefore, since by definition of $Y_{\lambda}^k$ 
all maximal chains start with
$\lambda$, we have our claim.
\end{proof}

We are interested in proving properties of $Y^k_{\lambda}$ when 
$\lambda$ is a rectangular partition, denoted by:

\begin{definition}
The principal order ideal of $Y^k$ generated by the 
partition $m^n$ will be denoted
\begin{equation}
L^k(m,n) = \left\{ \mu : \mu \preceq (m^n) \right\}
\, .
\end{equation}
\end{definition}
Note that $m\leq k$ since $Y^k$ contains only elements 
of $\mathcal P^k$.  Further,

\begin{remark} \label{remsauvemoi}
From Property~\ref{hookk}, $L^k(m,n)=L(m,n)$, for $n \leq k-m+1$. 
Therefore, all cases are covered when considering $n \geq k-m+1$.
That is, all $L^k(m,n)$ distinct from $L(m,n)$, plus the non distinct
case
$L^{n+m-1}(m,n)=L(m,n)$.
\end{remark}

We shall prove that $L^k(m,n)$ is a graded, self-dual, and 
distributive lattice. Further, we shall conjecture that its 
rank-generating function is unimodal in certain cases. 
To start, following from Proposition~\ref{rank},

\begin{corollary}
$L^k(m,n)$ is graded of rank $mn$.
\label{recrank}
\end{corollary}

It will develop that for each $k$, the principal order ideal 
generated by $(m^n)$ in $Y^k$ is isomorphic to an induced 
subposet of the principal order ideal generated by $(m^n)$ 
in the Young Lattice.  From this, we can then derive a number 
of properties for the posets $L^k(m,n)$.  To this end, we first 
explicitly determine the vertices of the order ideal.

\begin{theorem} 
\label{kvert}
The set of partitions in $L^k(m,n)$ are
those that fit inside an $m\times n$ rectangle and
have no more than $k-m+1$ rows shorter than $m$.
\end{theorem}

\begin{proof}
We first show that any $\lambda \in L^k(m,n)$ can have at most 
$k-m+1$ rows of length shorter than $m$.  
Suppose the contrary, and note that the top $k-m+2$ rows of 
$\lambda/^k$ form a partition since the rows are all
shorter than $m$ and thus no hook exceeds $k$.
Therefore, the first column of $\SS$ has height at least 
$k-m+2$ and $\lambda^{\omega_k}$ has a row of length at least 
$k-m+2$ by definition of $k$-conjugate.
However, the rows of $(m^n)^{\omega_k}$ do not exceed $k-m+1$
by Proposition~\ref{kconrec} and thus
$\lambda^{\omega_k}\not\subseteq (m^n)^{\omega_k}$.
Therefore $\lambda \not \in L^k(m,n)$ by Property~\ref{contains}.

On the other hand, to prove that any $\lambda=(m^a,\mu)\subseteq m^n$ 
with $\mu\subseteq (m^{k-m+1})$ lies in $L^k(m,n)$, it suffices to prove 
$\lambda\in L^k(m,\ell(\lambda))$ since $m^{\ell(\lambda)}\preceq m^{n}$
-- i.e. using Corollary~\ref{removebox} $m$ times,
we obtain $m^{n-1} \preceq m^{n}$, and by iteration
$m^{\ell(\lambda)} \preceq m^{n}$.  To this end, note that the top 
$\ell(\mu)\leq k-m+1$ rows of $(m^{\ell(\lambda)})/^k$ fit inside the 
shape $(m^{k-m+1})$ implying every diagonal has a distinct 
$k+1$-residue by Remark~\ref{recres}.  Therefore, Theorem~\ref{khook}
implies removable corners can be successively removed from
the top $\ell(\mu)$ rows in $(m^{\ell(\lambda)})$ to
obtain partitions $\lambda^{(1)},\ldots,\lambda^{(i)}$
where $\lambda \cpreceq \lambda^{(1)}  
\cpreceq \cdots \cpreceq \lambda^{(i)} \cpreceq (m^{\ell(\lambda)})$.
\end{proof}

The theorem reveals that when $n \geq k-m+1$,
the elements of $L^k(m,n)$ are of the form
$(m^a,\mu)$ for $\mu \subseteq (m^{k-m+1})$ and 
$a \leq n-(k-m+1)$.  By  Remark~\ref{remsauvemoi},
we have thus identified all the vertices.

\begin{corollary}
For $n \geq k-m+1$ and $m\leq k$, the set of partitions in $L^k(m,n)$ is
the disjoint union
\begin{eqnarray}
\left\{\mu\subseteq (m^{k-m+1})\right\}
\bigcup_{i=1}^{n-(k-m+1)}  
\left\{
\left(m^{i},\mu_1+1,\ldots,\mu_{k-m+1}+1\right)
: \mu\subseteq (m-1)^{k-m+1}
\right\} \nonumber
\, .
\end{eqnarray}
\label{kposet}
\end{corollary}

It also follows from Theorem~\ref{kvert} that the vertices of $L^m(m,n)$ 
are simply the partitions with at most one row smaller than $m$.

\begin{corollary}
The set of partitions in $L^m(m,n)$ is 
\begin{equation}
L^m(m,n) = 
\left\{(m^j,i) \,: 
\,
\text{$0\leq i\leq m$ and $0\leq j\leq n-1$}
\right\}
\,.
\label{justchain}
\end{equation}
\end{corollary}

\section{Further properties of $k$-Young lattice ideals}

Equipped with a simple characterization of the vertex set 
of $L^k(m,n)$, we can now investigate a connection between 
the $k$-Young lattice and the Young lattice.
As it turns out, the principal order ideal generated by $(m^n)$ 
in $Y^k$ is isomorphic to the induced subposet of $L(m,n)$
containing only the subset of partitions with no more than
$k-m+1$ rows smaller than $m$.  Consequently, 
it is easy to grasp which elements are covered by 
$\lambda$ in $L^k(m,n)$ and to deduce
that the posets are self-dual and distributive.

\begin{proposition}
\label{coverrela} Let $\lambda,\mu \in 
L^k(m,n)$.  Then $\mu\cpreceq\lambda$ if and only if 
$\mu <\!\!\!\cdot \, \lambda$.
\end{proposition}

\begin{proof}
First consider $\lambda\in L^k(m,n)$ where 
$\lambda\subseteq(m^{k-m+1})$.  Since $h_{(1,1)}(\lambda)\leq k$,
$\mu\preceq\lambda$ reduces to $\mu\leq\lambda$ by Property~\ref{kord}.
Any other element of $L^k(m,n)$ has the form 
$\lambda=(m^b,\nu)$ for some $\nu\in\mathcal P^m$ with 
$\ell(\nu)=k-m+1$ by Corollary~\ref{kposet}.  
Given $\lambda$ of this form, since
$\mu\cpreceq\lambda$ and $\mu<\!\!\cdot\lambda$
both require that $\mu=\lambda-e_r$ where
$r$ is a row of $\lambda$ with a removable corner,
we need only consider $b\leq r\leq b+\ell(\nu)$. 
Further, $(\lambda-e_r)\not\in L^k(m,n)$ if $r=b$ 
(the partition would have more than $k-m+1$ rows shorter
than $m$).   Thus it suffices to show that for 
$b<r\leq b+\ell(\nu)$, there is a removable corner in 
row $r$ of $\lambda$ if and only if
$(\lambda-e_r)\cpreceq\lambda$
-- equivalently by Theorem~\ref{khook} --
if and only if there is a removable corner in row 
$r$ of $\SS$ whose $k+1$-residue does not occur
in any higher removable corner.

When row $r$ of $\SS$ has a removable corner that is the highest 
of a given $k+1$-residue, there is a removable corner in row $r$ 
of $\lambda$ by Property~\ref{newcorner}.
On the other hand, if there is a removable corner in row 
$b<r\leq b+\ell(\nu)$ of $\lambda$ then
there is also a removable corner in this row of $\SS$ 
since the top $\ell(\nu)$ rows of $(m^b,\nu)/^k$ 
coincide with the diagram of $\nu$ given
that $\ell(\nu)=k-m+1$.
Furthermore, $\nu \subseteq (m^{k-m+1})$ also implies that
there cannot be another removable corner above the one
in row $r$ of the same $k+1$-residue because the diagonals
in $(m^{k-m+1})$ all have distinct $k+1$-residue.
\end{proof}

\begin{theorem}
\label{indk}
$L^k(m,n)$ is isomorphic to the induced subposet of $L(m,n)$
with vertices restricted to the partitions in $L^k(m,n)$.
Equivalently, for $ \lambda,\mu\in L^k(m,n)$, 
\begin{equation}
\mu\preceq\lambda\quad\iff\quad
\mu \subseteq\lambda \,.
\end{equation}
\label{recorderminus}
\end{theorem}

\begin{proof}
Let $\lambda,\mu \in L^k(m,n)$.  From Property~\ref{contains}, 
$\mu\preceq\lambda\implies\mu \subseteq \lambda$.  It thus
remains to show that there exists a chain from $\mu$ to 
$\lambda$ in $L^k(m,n)$ when $\mu \subseteq \lambda$. 
Equivalently by the previous proposition, it suffices to show 
that we can reach $\mu$ by successively adding boxes to $\lambda$ 
in such a way that no intermediate step gives a partition with 
more than $k-m+1$ rows of length less than $m$.  
This is achieved as follows:
given $\mu \subseteq \lambda$ in $L^k(m,n)$,
consider the chain of partitions $\mu=\mu^0\subseteq \mu^1
\subseteq\cdots\subseteq\mu^j=\lambda$ where
$\mu^{i+1}$ is obtained by adding one box to the
first row in $\mu^i$ that is strictly less than 
the corresponding row in $\lambda$.
Since the chain starts from $\mu=(m^a,\nu)$ where 
$\ell(\nu)\leq k-m+1$, the number of 
rows of length less than $m$ does not exceed $k-m+1$ in 
any $\mu^i$  by construction.
\end{proof}

We can now derive a number of properties for the order ideals $L^k(m,n)$
based on the identification with induced subposets of $L(m,n)$ under 
inclusion of diagrams.  First, given the explicit description 
Eq.~\eqref{justchain} for the vertices in $L^m(m,n)$ we have

\begin{proposition}
\label{justchainposet}
The order ideal $L^m(m,n)$ is isomorphic to
the saturated chain of partitions:
$$
\emptyset\subseteq (1) \subseteq 
\cdots\subseteq (m)\subseteq (m,1)\subseteq\cdots
\subseteq
(m^j,i) \subseteq\cdots
\subseteq (m^{n-1},m-1) \subseteq (m^n)
\,.
$$
\end{proposition}

\begin{proposition}
\label{inducedposet}
For $k\geq m$, $L^{k}(m,n)$ is an induced subposet of $L^{k+1}(m,n)$.
\end{proposition}
\begin{proof}
From Theorem~\ref{kvert}, the elements of $L^{k}(m,n)$ (or $L^{k+1}(m,n)$)
are partitions contained in $(m^n)$ with at most $k-m+1$ (resp. $k-m+2$) rows 
that are smaller than $m$.  
Therefore, by Theorem~\ref{recorderminus}, it suffices to note that
under inclusion of diagrams, the poset of partitions 
that fit inside an $m\times n$ rectangle with no more 
than $k-m+1$ parts smaller than $m$ is an induced
subposet of the poset of partitions 
that fit inside an $m\times n$ rectangle with no more 
than $k-m+2$ parts smaller than $m$.
\end{proof}

\begin{proposition}  Let $\lambda,\mu\in L^k(m,n)$.

\noindent (i) $\lambda\wedge\mu$ is the partition determined by
the intersection of the cells in the diagrams in $\lambda$ and $\mu$.

\noindent (ii) $\lambda\vee\mu$ is the partition whose diagram 
is determined by the union of the cells in $\lambda$ and $\mu$.
\label{projoin}
\end{proposition}

\begin{proof}  
Since the meet and join of elements in $L(m,n)$
is given by the intersection and union of diagrams respectively, 
and $L^k(m,n)$ is isomorphic to an induced subposet of $L(m,n)$ by
Theorem~\ref{indk}, it suffices to show that $L^k(m,n)$ is
closed under the intersection and the union of diagrams.
Equivalently, by Theorem~\ref{kvert}, we must prove that the 
intersection and union of such partitions do not have 
more than $k-m+1$ rows with length less than $m$.
Let $\lambda=(m^a,\bar \lambda)$ with $\ell(\bar \lambda) \leq k-m+1$ 
and $\mu=(m^b,\bar \mu)$ with $\ell(\bar \mu) \leq k-m+1$
and assume $a\geq b$ without loss of generality. 
As such, the diagram of $\lambda$ intersected with $\mu$ 
has at least $b$ rows of length $m$, and at most 
$\ell(\bar \mu)$ rows with length less than $m$.
Therefore, there are no more than $k-m+1$ rows of length 
smaller than $m$.  Similarly, the union has at least $a$ rows of 
length $m$, and at most 
$\max\bigl \{\ell(\bar\lambda) ,\ell(\bar\mu)-(a-b) \bigr \}$  
rows less than $m$.
Again, no more than $k-m+1$ rows of length smaller than $m$.
\end{proof}

Now given that the meet and join  of $L^k(m,n)$ coincide with 
those of $L(m,n)$.  Therefore, since $L^k(m,n)$ is an induced  
subposet of the {\it lattice} $L(m,n)$, we have\footnote{Note that the
corollary also
follows immediately from the fact that $Y^k$ is a lattice.}

\begin{corollary} 
For each $k\geq m$, $L^k(m,n)$ is a lattice.
\end{corollary} 

Furthermore, since $L(m,n)$ is a distributive lattice,
the relations 
\begin{equation}
\lambda \vee (\mu \wedge \nu) = 
(\lambda \vee \mu) \wedge (\lambda \vee \nu) \, ,
\qquad
\lambda \wedge (\mu \vee \nu) = 
(\lambda \wedge \mu) \vee (\lambda \wedge \nu) \, , 
\end{equation}
must hold.  Therefore,  these relations hold in 
the induced subposets $L^k(m,n)$ and we find

\begin{corollary}
For each $k\geq m$, $L^k(m,n)$ is distributive.
\end{corollary}

In addition to having that each of the $L^k(m,n)$ are
distributive lattices, we can also prove that they are
symmetric.

\begin{theorem}
For any $k\geq m$, $L^k(m,n)$ is self-dual.
\end{theorem}
\begin{proof}
Let $\bar L^k(m,n)$ denote the dual of $L^k(m,n)$ 
and consider the mapping $\phi(\lambda)=\lambda^c$
where $\lambda^c$ is the partition determined by rotating
the complement of $\lambda$ in $(m^n)$ by $180^{o}$.
Since $\lambda\in L^k(m,n)$ implies that $\lambda=(m^a,\mu)$ for some 
$a\leq n-(k-m+1)$ and $\mu\subseteq (m^{k-m+1})$ by Corollary~\ref{kposet},
$\lambda^c=(m^{n-(k-m+1)-a},\mu^c)$ for some 
$\mu^c\subseteq (m^{k-m+1})$ satisfying $|\lambda|+|\lambda^c|=mn$.  
Therefore $\phi: L^k(m,n)\to\bar L^k(m,n)$ and it suffices to 
show that $\phi$ is an order-preserving bijection.  Equivalently,
that $\mu\preceq\lambda$ in $L^k(m,n)
\iff \phi(\lambda)\preceq\phi(\mu)$ in $\bar L^k(m,n)$.
By Theorem~\ref{indk} and the definition of $\bar L^k(m,n)$, 
this is equivalent to $\mu \subseteq\lambda \iff \lambda^c\subseteq\mu^c$
for elements $\lambda,\mu \in L^k(m,n)$, which is true.
\end{proof}

\begin{corollary}
$L^k(m,n)$ is rank-symmetric.
\end{corollary}

\section{Rank-generating function}

The explicit description of the partitions in the posets $L^k(m,n)$ 
can also be used to determine the rank-generating functions.  
Recall that the number of elements of rank $i$ in $L(m,n)$ is 
the coefficient of $q^i$ in the Gaussian polynomial:
\begin{equation} 
\sum_{i=0}^{mn} p_i(m,n)\,q^i = 
\begin{bmatrix}
m+n\\m
\end{bmatrix}_q
\, .
\label{binid}
\end{equation}
Similarly, we can determine the rank-generating functions for
$L^k(m,n)$.

\begin{theorem}
For $n \geq k-m+1$ and $m\leq k$, the number of elements of degree $i$ 
in $L^{k}(m,n)$ is the coefficient of $q^i$ in
\begin{equation}
\label{pk}
\sum_{i=0}^{mn} p^{k}_i(m,n)\,q^i
\,:=\,
\begin{bmatrix}
k+1 \\m
\end{bmatrix}_q
+
q^{k+1} \frac{1-q^{m(n-k+m-1)}}{1-q^m}
\begin{bmatrix}
k \\m-1
\end{bmatrix}_q \,. 
\end{equation}
\end{theorem}
\begin{proof}
Recall that Corollary~\ref{kposet} provides an
interpretation for the vertices of $L^k(m,n)$ 
as a disjoint union of sets whose elements can
be understood in terms of certain $L(a,b)$:
\begin{eqnarray}
\left\{\mu\subseteq (m^{k-m+1})\right\}
\bigcup_{i=1}^{n-(k-m+1)}  
\left\{
\left(m^{i},\mu_1+1,\ldots,\mu_{k-m+1}+1\right)
: \mu\subseteq (m-1)^{k-m+1}
\right\} \nonumber
\, .
\end{eqnarray}
Identity~\eqref{binid} then gives that
the number of elements of rank $i$ in $L^k(m,n)$
is the coefficient of $q^i$ in
\begin{equation}
\begin{bmatrix}
k+1 \\m
\end{bmatrix}_q
+
\sum_{i=1}^{n-(k-m+1)}
q^{m(i)+k-m+1}
\begin{bmatrix}
k \\m-1
\end{bmatrix}_q
\, .
\label{easier}
\end{equation}
\end{proof}

Letting $q\to 1$ in Eq.~\eqref{easier}
then gives the number of vertices:

\begin{corollary}
For $n \geq  k-m+1$ and $m\leq k$, the number of vertices in $L^k(m,n)$ is
\begin{equation}
\sum_{i= 0}^{mn} p^k_i(m,n) 
\,=\,
\left(
\begin{matrix}
k+1 \\m
\end{matrix}
\right)
+\left( n-k+m-1\right)
\begin{pmatrix}
k \\m-1
\end{pmatrix}
\, .
\end{equation}
\end{corollary}

In the next section, we will see that the study of the order ideals 
$L^k(m,n)$ is simplified by a close examination of certain
subsets of the vertices that are determined by $k$.

\begin{definition}
For $n \geq k-m+1$ and $k>m$, let $\Gamma^k(m,n)$ denote
the set of all partitions in $L^{k}(m,n)$ that are not in 
$L^{k-1}(m,n)$.
\end{definition}

We conclude the section with a discussion of the sets
$\Gamma^k(m,n)$, starting with an explicit description 
of the elements.

\begin{proposition}
The elements of $\Gamma^k(m,n)$
are partitions that fit inside an $m\times n$
rectangle and have exactly $k-m+1$ rows
with length smaller than $m$.  Equivalently, 
\label{krecur} 
\begin{equation}  \label{quotient}
\Gamma^k(m,n)
=
\Bigl\{(m^a,\mu_1+1,\dots,\mu_{k-m+1}+1) \,: \, 
0\leq a\leq n-(k-m+1) \, , \, 
\mu \subseteq (m-2)^{k-m+1}
\Bigr\} 
\end{equation}
\end{proposition}
\begin{proof}
Theorem~\ref{kposet} indicates that the vertices of 
$L^k(m,n)$ are partitions with no more than $k-m+1$ rows of length 
smaller than $m$ while those of $L^{k-1}(m,n)$ are partitions
with no more than $k-m$ rows of length smaller than $m$.
The result thus follows.
\end{proof}

Since the vertices of $L^m(m,n)$ are given in Eq.~\eqref{justchain} to be 
the partitions with no more than one row of length smaller than $m$,
the set of vertices for  $L^k(m,n)$ can be obtained by adding the 
elements of $\Gamma^{i}(m,n)$ for $i=m+1,\ldots,k$
to this set of partitions.  That is,

\begin{proposition}
For $n\geq k-m+1$ and $k > m$, the vertices of $L^k(m,n)$ 
are the disjoint unions:
\begin{equation}
L^k(m,n) = L^m(m,n)\;\cup_{j=m+1}^k\Gamma^j(m,n)
\, .
\end{equation}
\label{beginrec}
\end{proposition}

\begin{figure}[htb]
\begin{center}
\epsfig{file=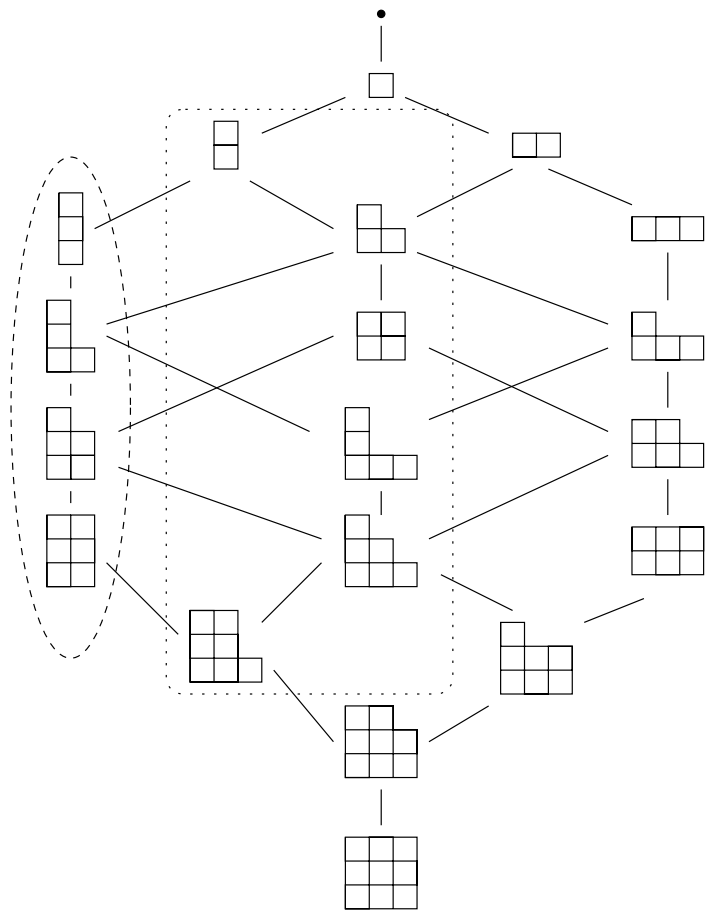}
\end{center}
\caption{Decomposing $L^5(3,3)$ into 
$L^3(3,3)\cup\Gamma^4(3,3)\cup\Gamma^5(3,3)$}
\label{trois}
\end{figure}

Given this decomposition for the set of vertices in $L^k(m,n)$,
we note that there are no edges between partitions in $\Gamma^j(m,n)$ 
and partitions in $\Gamma^{j-i}(m,n)$ for any $i>1$ since adding or 
deleting a box to a partition with exactly $a$ rows smaller than 
length $m$ never gives rise to a partition with more than $a+1$ 
rows smaller than $m$.  Now, $L^{m+1}(m,n)$ can be constructed 
by connecting elements of $\Gamma^{m+1}(m,n)$ to the 
appropriate vertices in the saturated chain $L^m(m,n)$ of 
Proposition~\ref{justchainposet}.  $L^{m+2}(m,n)$ can then be 
constructed by connecting {\it only} elements of $\Gamma^{m+1}(m,n)$ 
to the appropriate elements of $\Gamma^{m+2}(m,n)$.  
In this manner, the poset $L^k(m,n)$ can be 
constructed from the saturated chain $L^m(m,n)$.  In particular, 
$L(m,n)= L^{m+n-1}(m,n)$ can be obtained using this process
(see Figure~\ref{trois} for an example).  In summary,

\begin{remark}
The poset $L^k(m,n)$ can be obtained from $L^{k-1}(m,n)$ by adding 
edges from $\lambda\in\Gamma^{k-1}(m,n)$ 
to all partitions in $\Gamma^{k}(m,n)$ 
that contain or are contained in $\lambda$.
\end{remark}

We now obtain the rank-generating function for 
$\Gamma^k(m,n)$ from Proposition~\ref{krecur} by using
the Gaussian polynomial for partitions inside an
$(k-m+1)\times (m-2)$ rectangle.

\begin{proposition}  \label{ucoro}
For  $n\geq k-m+1$  and $k>m$, the rank-generating function
for $\Gamma^k(m,n)$ is
\begin{equation} 
\label{useq}  \sum_{i \geq 0} u_i^k(m,n) \,q^i
\,
:= \, \sum_{\lambda\in\Gamma^k(m,n)} q^{|\lambda|}
\, = \, q^{k-m+1} 
\frac{(1-q^{m(n-k+m)})}{(1-q^m)} 
\begin{bmatrix}
k-1 \\m-2
\end{bmatrix}_q \, .
\end{equation}
In this equation,  $u_i^k(m,n)$ is defined to be 
the number of elements 
of degree $i$ in $\Gamma^{k}(m,n)$.
\end{proposition}

It turns out that the sequence of coefficients in this expression
is rank-symmetric.

\begin{proposition} 
\label{ranksym}
$u^k_i(m,n)=u^k_{mn-i}(m,n)$ for all 
$i=0,\ldots,\lfloor \frac{mn}{2} \rfloor$.
That is, the vector of coefficients 
$\vec u^k=(u_0^k(m,n),\ldots,u_{mn}^k(m,n))$ 
is symmetric about the middle ({\it i.e.} $mn/2$).
\end{proposition}
\begin{proof}
By Eq.~\eqref{useq},
$u^k_i(m,n)$ is the coefficient of $q^i$ in a product of
three rank-symmetric polynomials:
$q^{k-m+1}$, $\frac{(1-q^{m(n-k+m)})}{(1-q^m)}$, and
$\begin{bmatrix}
k-1 \\m-2
\end{bmatrix}_q$.
Therefore, since the term of lowest degree in the expansion 
of these three polynomials has degree $k-m+1$ and the highest 
degree term has degree $k-m+1 + m(n-k+m-1)+ (m-2)(k-m+1)=mn-k+m-1$,
the polynomial must be rank-symmetric about 
$\bigl(k-m+1 +(mn-k+m-1)\bigr)/2=mn/2$.
\end{proof}

\section{Unimodality and sieved sums}
\label{uni}

Recall that a poset $P$ of rank $d$ is {\bf unimodal} if
$$
p_0\leq p_1\leq\cdots\leq p_i\geq p_{i+1}\geq\cdots\geq p_{d}
\text{ for some } 0\leq i\leq d\,,
$$
where $p_i$ is the number of elements with rank $i$ in $P$.
We shall say the rank-generating function for $P$ is unimodal when 
the vector $(p_0,\ldots,p_d)$, with $p_i$ the 
coefficient of $q^i$, forms a unimodal sequence.  
For example, the Gaussian polynomial \eqref{binid} is unimodal 
because it is known \cite{[P],[S]} that the coefficients 
$(p_0,\ldots,p_{mn})$ form a unimodal sequence.
Equivalently, $L(m,n)$ is  a unimodal poset -- one of the deeper 
properties of this order ideal. 

In this section, we address the question of unimodality 
for the posets $L^k(m,n)$, illustrating the role of $k$ in 
our study.  As such, by Eq.~\eqref{pk}, the order ideal 
$L^k(m,n)$ will be unimodal when
\begin{equation}
\sum_{i=0}^{mn}
p_i^{k}(m,n)\, q^i
\,=\,
\begin{bmatrix}
k+1 \\m
\end{bmatrix}_q
+
q^{k+1} \frac{1-q^{m(n+m-k-1)}}{1-q^m}
\begin{bmatrix}
k \\m-1
\end{bmatrix}_q \,
\label{kbin}
\end{equation}
is unimodal.  
Our work in the preceding section pays off here by
enabling us to rewrite this expression as a sum
of rank-symmetric components.

\begin{proposition} 
\label{prop1plus}
The number of elements in $L^k(m,n)$ is given by
the coefficient of $q^i$ in
\begin{equation}
\sum_{i=0}^{mn} p_i^{k}(m,n)\,q^i 
=  \sum_{i=0}^{mn}
\left( 1 \; + \; \sum_{r=m+1}^{k}
u_i^r(m,n)\right) q^i \, ,
\label{1plus}
\end{equation}
where $u_i^r(m,n)$ is defined in Eq.~\eqref{useq}.
\end{proposition}
\begin{proof}
Proposition~\ref{beginrec} decomposes the set of 
vertices of $L^k(m,n)$ into $L^m(m,n)\cup_{i=m+1}^k \Gamma^i(m,n)$,
implying that the coefficients $p_i^{k}(m,n)$ occur in the
expression:
\begin{equation}
\sum_{i=0}^{mn} p_i^{k}(m,n)\,q^i 
 = 
\sum_{\lambda\in L^m(m,n)} q^{|\lambda|}
+\sum_{r=m+1}^k \sum_{\lambda\in\Gamma^r(m,n)}
q^{|\lambda|} \, . 
\end{equation}
The result then follows since $L^m(m,n)$ is just a saturated 
chain by Proposition~\ref{justchainposet}, and the elements 
in $\Gamma^r(m,n)$ are given in Proposition~\ref{ucoro}
by $u_i^r(m,n)$.
\end{proof}

Since $u_i^r(m,n)=u_{mn-i}^r(m,n)$ by Proposition~\ref{ranksym},
\eqref{1plus} reveals 
$\vec p^{\, k}(m,n) =(p_0^k(m,n),\ldots,p_{mn}^k(m,n))$
is a sum of sequences symmetric about the middle.
Thus, the question of unimodality of $L^k(m,n)$ reduces to 
a study of the sequences $\vec u^{\,r}(m,n)$.
We have experimentally discovered that certain sums of 
these sequences are unimodal under conditions 
on $m$ expressed in Conjecture~\ref{conjecu}, and
more generally in Conjecture~\ref{genconj}.
This given, we can deduce that $\vec p^{\, k}(m,n)$
is unimodal in these cases and can be written as a sum of 
rank-unimodal sequences (see Remark~\ref{remunisum}).
To this end, we start with the special case when $m$ is prime.

\begin{conjecture} \label{conjecu}
If $n \geq k-m+1$ and $k \neq -1,0 \mod m$ for a 
prime number $m<k$, then
\begin{equation}
\sum_{i=0}^{mn} u_i^{k}(m,n)
\, q^i
\,=\,
q^{k-m+1}\,
\frac{(1-q^{m(n-k+m)})}{(1-q^m)} 
\begin{bmatrix}
k-1 \\m-2
\end{bmatrix}_q
\end{equation}
is unimodal. Further, when $n > k-m+1$ and $k = -1 \mod m$,
\begin{eqnarray} 
& &\sum_{i=0}^{mn}\Bigl(u_i^{k}(m,n)+u_i^{k+1}(m,n)\Bigr) q^i 
  \nonumber \\
& & \qquad \,=\, 
q^{k-m+1}\,
\frac{(1-q^{m(n-k+m)})}{(1-q^m)} 
\begin{bmatrix}
k-1 \\m-2
\end{bmatrix}_q+
q^{k-m+2}\,
\frac{(1-q^{m(n-k+m-1)})}{(1-q^m)} 
\begin{bmatrix}
k \\m-2
\end{bmatrix}_q
\end{eqnarray}
is unimodal. 
\end{conjecture}

There are two consequences of this conjecture;
the first relating to the unimodality of $L^k(m,n)$.

\noindent
{\bf Consequence of Conjecture~\ref{conjecu}.}
{\it
\label{conjecp}
If $n>k-m+1$ and $k \neq -1 \mod m$ for a 
prime $m<k$, then the order ideal $L^k(m,n)$ is unimodal.
Equivalently,
\begin{equation}
\sum_{i=0}^{mn} p_i^{k}(m,n)
\, q^i
\,=\,\begin{bmatrix}
k+1 \\m
\end{bmatrix}_q
+
q^{k+1} \frac{1-q^{m(n+m-k-1)}}{1-q^m}
\begin{bmatrix}
k \\m-1
\end{bmatrix}_q \,
\end{equation}
is unimodal under these conditions.
}

\smallskip

\noindent
{\it Proof following from Conjecture~\ref{conjecu}.} 
By Eq.~\eqref{1plus}, it suffices to show that
$\sum_{j=m+1}^k \vec u^j(m,n)$ is unimodal.
For all $m<j<n+m-1$, the sequences $\vec u^j(m,n)$ are rank-symmetric 
about $mn/2$ by Proposition~\ref{ranksym}.  Therefore, 
if $\vec u^j(m,n)$ is unimodal for $j \neq -1,0 \mod m$ and
$\vec u^j(m,n)+{\vec u}^{j+1}(m,n)$ is unimodal when $j=-1 \mod m$, 
then $\sum_{j=m+1}^k \vec u^j(m,n)$ is unimodal unless $k=-1 \mod m$.
That is, except in the case that ${\vec u}^{j+1}(m,n)$
cannot be added to restore unimodality.
\hfill$\square$

Interestingly, a second consequence of Conjecture~\ref{conjecu} 
relates to sieved $q$-binomial coefficients (see for example,
\cite{[GS],[SW],[WW]}).   In this result, 
``the sum of the coefficients of $q^{\ell+*m}$'' in an
expression refers to the sum of coefficients of $q^i$ 
for every $i=\ell \mod m$.

\begin{proposition} \label{proppart}
If $k \neq -1,0 \mod m$ for a prime $m<k$,
then for any $0\leq\ell\leq m-1$,
the sum of the coefficients of $q^{\ell+*m}$ in 
$\begin{bmatrix} k-1\\m-2\end{bmatrix}_q$
is $\begin{pmatrix}k-1\\m-2\end{pmatrix}/m$.
Equivalently,
\begin{equation}
\sum_{j\geq 0}p_{\ell+jm}(m-2,k-m+1)= 
\frac{1}{m}\begin{pmatrix}k-1\\m-2\end{pmatrix}\,, \quad\text{for fixed }\;
\ell=0,\dots,m-1\, .
\end{equation}
\end{proposition}

We shall demonstrate how this proposition follows from 
Conjecture~\ref{conjecu} to give evidence supporting 
the conditions under which we believe unimodality to hold.
However, the result is implied from the main result in 
\cite{[WW]} (and more recently from \cite{[SW]})
as follows:
\begin{proof}
Condition 3 of the main result in \cite{[WW]} says that
the sum of the coefficients of $q^{\ell+*m}$ in 
$\begin{bmatrix} n\\t \end{bmatrix}_q$
is $\begin{pmatrix}n\\t \end{pmatrix}/m$ for all $\ell$ iff 
for all $d>1$ that divide $m$, we have $t \mod d > n \mod d$.
Since only $d=m$ divides $m$ when $m$ is prime, 
letting $n=k-1$ and $t=m-2$, we note that 
$m-2 \mod m > k-1 \mod m$ is equivalent to 
the condition $k \neq -1,0 \mod m$.
\end{proof}

\noindent
{\it Proof following from Conjecture~\ref{conjecu}.} 
For $k \neq -1,0 \mod m$ for $m$ prime, 
Conjecture~\ref{conjecu} implies
$$
\frac{(1-q^{m(n-k+m)})}{(1-q^m)} 
\begin{bmatrix}
k-1 \\m-2
\end{bmatrix}_q
$$
is unimodal.  Letting 
$n \to \infty$ in this expression then produces an 
increasing (unimodal and symmetric about infinity)
sequence $(v_0^k(m),v_1^k(m),\ldots)$, where
\begin{equation} 
\sum_i v_i^k(m) \, q^i = 
(1+q^m+q^{2m}+\cdots) 
\begin{bmatrix}
k-1 \\m-2
\end{bmatrix}_q 
\, .
\end{equation}
Using the definition of $p_i(m,n)$ in Eq. \eqref{binid},
the coefficient of $q^i$ in this expression satisfy
\begin{equation*}
v_i^k(m) = p_{i}(m-2,k-m+1) + p_{i-m}(m-2,k-m+1) +p_{i-2m}(m-2,k-m+1) 
+ \cdots \, .
\end{equation*}
This sequence thus meets the conditions in the following proposition with
$\vec a=\vec p(m-2,k-m+1)$, where the claim follows from Eq.~\eqref{ik}
since $\sum_i p_i(m-2,k-m+1)=\begin{pmatrix} k-1\\m-2 \end{pmatrix}$.
\hfill $\square$

\begin{proposition} \label{propmax}
Let $v_i = \sum_{t\geq 0} a_{i-t m}$
for some $\vec a=(a_0,\ldots,a_D)$ and
$a_j=0$ when $j<0$ or $j>D$. 
If the infinite sequence $\vec v=(v_0,v_1,\ldots,)$
is weakly increasing, then

\noindent (i) $v_i=v_{D-m+1}$ 
for $i \geq D-m+1$.

\noindent (ii) 
For any $\ell=0,\dots,m-1 $,
\begin{equation}
\label{ik}
\sum_{j\geq 0} a_{\ell+jm} = 
\frac{1}{m}\sum_{j=0}^{D} a_j
\, .
\end{equation}
\end{proposition}

\begin{proof} 
\noindent (i) 
Since $i+m>D$ when $i \geq D-m+1$, 
we have $a_{i+m}=0$.
The definition of $v_{i+m}$ implies
$v_{i+m}=a_{i+m}+v_i$ and thus $v_{i+m}=v_i$ when $i \geq D-m+1$.
To show $v_i=v_{D-m+1}$ for all $i\geq D-m+1$,
note that
$v_{D-m+1} \leq v_{D-m+2} \leq \cdots\leq v_{D+1}$
since $\vec v$ is increasing.  However, $v_{D-m+1}=v_{D+1}$ by 
the preceding argument, implying $v_{D-m+1} =v_{D-m+2}\cdots=v_{D+1}$.  
By iteration, $v_i=v_{D-m+1}$ for all $i \geq D-m+1$.

\noindent (ii) 
Since $a_{i}=0$ for $i<0$, we have for all $1\leq \ell\leq  m$ that
\begin{equation}
v_{D-\ell+1} = a_{D-\ell+1} + a_{D-\ell+1-m}+ \cdots + 
a_{D-\ell+1-\lfloor (D-\ell+1)/m \rfloor  m} \, .  
\label{righthere}
\end{equation}
That is, 
$v_{D-\ell+1}=\sum_{n} a_n$ for $n=D-\ell+1 \mod m$.
The sum $v_{D-m+1}+\dots+v_{D-1}+v_D$ is thus the
sum of all entries of $\vec a$.  
However, (i) implies that $v_{D-m+1}=\cdots=v_{D-1}=v_D$ 
and therefore, $m\cdot v_{D-\ell+1}=\sum_{j=0}^D a_j $,
for $\ell=1,\dots,m$.  From Eq.~\eqref{righthere}, $v_{D-\ell+1}$
is the sum of entries in $\vec a$ indexed by $D-\ell+1$ modulo $m$,
thus implying the claim since as $\ell$ runs over 1 to $m$, 
$D-\ell+1$ runs over all possible values modulo $m$.
\end{proof}

In summary, when $m$ is prime, proving unimodality of the posets 
$L^k(m,n)$ reduces to proving the unimodality of
$\vec u^{\, r}(m,n)+\vec u^{\, r+1}(m,n)$.  As such,
to extend this idea for $m$ not prime, we study
the unimodality of more general sums of $\vec u^{\,r}(m,n)$.

\begin{conjecture} \label{genconj}
Consider integers $a$ and $b$ such that $m\leq a<b\leq n+m-1$.
If $a,b \neq -1 \mod p$ for every prime divisor $p$ of $m$,
then 
\begin{equation}
\label{genconjeq}
\sum_{i=0}^{mn}
\left(
\sum_{j=a+1}^b 
u_i^j(m,n)
\right)q^i 
\,=\,
\sum_{j=a+1}^b q^{j-m+1} 
\frac{(1-q^{m(n-j+m)})}{(1-q^m)} 
\begin{bmatrix}
j-1 \\m-2
\end{bmatrix}_q \, 
\end{equation}
is unimodal.
\end{conjecture}

We recover Conjecture~\ref{conjecu} by taking $m$ prime
and letting $a=k-1$ and $b=k$ to get the unimodality of
$\vec u^k(m,n)$. Further, $a=k-1$ and $b=k+1$ implies
the unimodality of $\vec u^k(m,n)+\vec u^{k+1}(m,n)$.
As with the case of $m$ prime, we can extract two consequences
of Conjecture~\ref{genconj}.

\noindent
{\bf Consequence of Conjecture~\ref{genconj}.} 
{\it Let $k>m$. 
If $k \neq -1 \mod p$ for every prime divisor $p$ of 
$m$, 
then $L^k(m,n)$ is unimodal for all $n > k-m+1$. Equivalently,
\begin{equation}
\begin{bmatrix}
k+1 \\m
\end{bmatrix}_q
+
q^{k+1} \frac{1-q^{m(n+m-k-1)}}{1-q^m}
\begin{bmatrix}
k \\m-1
\end{bmatrix}_q \,
\end{equation}
is unimodal under these conditions.
}

\noindent
{\it Proof following from Conjecture~\ref{genconj}.} 
Given $k \neq -1 \mod p$ for all prime divisors $p$ of $m$
and noting that $m \neq -1 \mod p$, we meet the conditions of 
Conjecture~\ref{genconj} with $a=m$ and $b=k$.  Therefore, 
$\vec p^{\, k}(m,n)$ is unimodal using
the decomposition given in Proposition~\ref{prop1plus}.
\hfill$\square$

\begin{remark} 
\label{remunisum}
Following trivially from Proposition~\ref{prop1plus}, we can 
obtain an alternative decomposition for $\vec p^{\,k}(m,n)$ in
terms of sums of $\vec u^{\,r}(m,n)$:
for nonnegative integers $k_1<k_2<\cdots< k_\ell=k$,
when each $k_j \neq -1\mod p$ for all prime divisors $p$ 
of $m\leq k$,
\begin{equation}
\sum_{i=0}^{mn} p_i^{k}(m,n)\,q^i 
=  \sum_{i=0}^{mn}
\left( 1 \; + \sum_{r={m+1}}^{k_1} u_i^r(m,n) +
 \sum_{r={k_1+1}}^{k_2} u_i^r(m,n)+
\cdots + \; \sum_{r=k_{\ell-1}+1}^{k_\ell}
u_i^r(m,n)\right) q^i \, .
\end{equation}
Thus, by Conjecture~\ref{genconj}, the order ideals $L^k(m,n)$
can be decomposed in terms of $\ell+1$ unimodal pieces
rank-symmetric about the same point $mn/2$.
\end{remark}

Our second consequence is a more general result on sieved 
$q$-binomial coefficients.  We state this consequence as 
a proposition since we will provide a proof that is
independent of our conjecture.

\begin{proposition}  \label{conj4}
If $m\leq a<b$ are non-negative integers where
$a,b \neq -1\mod p$ for every prime divisor $p$ of $m$,
then 
the sum of the coefficients of $q^{\ell+*m}$
in
$\sum_{j=a+1}^b q^{j-(a+1)} 
\begin{bmatrix}
j-1 \\m-2
\end{bmatrix}_q$
is 
$\frac{1}{m} \sum_{j=a+1}^{b} 
\begin{pmatrix}
j-1 \\m-2
\end{pmatrix}$,  for any $0\leq\ell\leq m-1$.  
\end{proposition}

Note that even when $m$ is prime, this 
is more general than Proposition~\ref{proppart}.  
For example, the case that $k=-1 \mod m$ is included
when $a=k-1$ and $b=k+1$.  Before proving this result,
we shall show how it follows from Conjecture~\ref{genconj} 
to support the validity of our conjecture.

\noindent
{\it Proof following from Conjecture~\ref{genconj}.}
Letting $n\to\infty$ in Eq.~\eqref{genconjeq},
(and multiplying by $q^{m-a-2}$ for convenience),
the coefficients $(v_0^{a,b}(m),v_1^{a,b}(m),\dots)$ 
defined by
\begin{equation}
\sum_{i\geq 0} v_i^{a,b}(m) \, q^i\,=\,
(1+q^m+q^{2m}+ \cdots)\sum_{j=a+1}^b q^{j-(a+1)} 
\begin{bmatrix}
j-1 \\m-2
\end{bmatrix}_q \, 
\end{equation}
form an increasing sequence given Conjecture~\ref{genconj}.
Thus, defining 
\begin{equation} \label{piab}
\sum_{i\geq 0} p_i^{a,b}(m) \, q^i=
\sum_{j=a+1}^b q^{j-(a+1)} 
\begin{bmatrix}
j-1 \\m-2
\end{bmatrix}_q \, ,
\end{equation}
we have
\begin{equation}
v_i^{a,b}(m) = p_i^{a,b}(m)+p_{i-m}^{a,b}(m)+p_{i-2m}^{a,b}(m) + \cdots
\, .
\end{equation}
The proposition thus follows from Proposition~\ref{propmax} with
$\vec a = (p^{a,b}_0(m),p^{a,b}_{1}(m),\ldots)$, since
$\sum_i p_i^{a,b}(m)=\sum_{j=a+1}^b \begin{pmatrix}
j-1 \\m-2
\end{pmatrix}$
.
\hfill$\square$

\begin{proof}
From 
\begin{equation}
\sum_{\omega: \omega^m=1} \omega^r = 
\begin{cases}
m & {\rm if~} m|r \\
0 & {\rm otherwise} 
\end{cases} \, ,
\end{equation}
we have that the sum of coefficients 
of $q^{\ell+*m}$ in any polynomial $P(q)$ is
\begin{equation}
\frac{1}{m} \sum_{\omega: \omega^m=1} \omega^{-\ell} P(\omega) \, .
\end{equation}
Thus, to prove the sum of the coefficients of $q^{\ell+*m}$ in 
$\label{modm}
\sum_{j=a+1}^b q^{j-(a+1)} 
\begin{bmatrix}
j-1 \\m-2
\end{bmatrix}_q$ is 
$\sum_{j=a+1}^{b} 
\begin{pmatrix}
j-1 \\m-2
\end{pmatrix}/m$,
it suffices to prove
\begin{equation}
\frac{1}{m} \sum_{j=a+1}^b \sum_{\omega: \omega^m=1}
 \omega^{-(\ell-j+a+1)} \begin{bmatrix}
j-1 \\m-2
\end{bmatrix}_{\omega} =\frac{1}{m} \sum_{j=a+1}^{b} 
\begin{pmatrix}
j-1 \\m-2
\end{pmatrix} \, .
\end{equation}
Or equivalently, 
since the right hand side equals the $\omega=1$ term 
in the left hand side,
\begin{equation}
\frac{1}{m} \sum_{j=a+1}^b \sum_{ \begin{subarray}{c}
\omega^m=1 \\ \omega \neq 1 \end{subarray}}
 \omega^{-(\ell-j+a+1)} \begin{bmatrix}
j-1 \\m-2
\end{bmatrix}_{\omega} = 0
\, .
\end{equation}
To this end, we shall demonstrate that
for all $\omega$  such that $\omega^m=1$ and $\omega \neq 1$,
\begin{equation} \label{ademont}
\sum_{j=a+1}^b 
\omega^{j} \begin{bmatrix}
j-1 \\m-2
\end{bmatrix}_{\omega} = 0
\end{equation}
by proving this identity holds
when $\omega$ is a primitive $d^{th}$ root of unity, for all
$d|m$ not equal to 1.

If $j\neq -1,0\mod d$, then
$\begin{bmatrix}
j-1 \\m-2
\end{bmatrix}_{\omega}=0$
since the numerator of
$\begin{bmatrix}
j-1 \\m-2
\end{bmatrix}_{\omega}$ 
has one more zero than its denominator, given that
$i\mod d=0$ for some $j-d+2\leq i\leq j-1$.
On the other hand, when $j = -1 \mod d$
(see Lemma~1(3) of \cite{[SW]})
\begin{equation}
\begin{bmatrix}
j-1 \\m-2
\end{bmatrix}_{\omega}=\begin{pmatrix}
n-1 \\m/d -1
\end{pmatrix} \, ,
\end{equation}
for $n=(j+1)/d-1$.  Thus, using
$$
\begin{bmatrix}
j \\m-2
\end{bmatrix}_{\omega} = \frac{1-\omega^{j}}{1-\omega^{j-m+2}}
\begin{bmatrix}
j-1 \\m-2
\end{bmatrix}_{\omega} \, ,
$$
we also have
\begin{equation}
\begin{bmatrix}
j \\m-2
\end{bmatrix}_{\omega}=  \frac{1-\omega^{-1}}{1-\omega}
\begin{pmatrix}
 n-1 \\m/d-1 
\end{pmatrix} \, .
\end{equation}
However, the conditions $a,b \neq -1 \mod p$ for every prime
divisor $p$ of $m$ imply that $a,b \neq -1 \mod d$ for any 
$d|m$.  Therefore, the cases $a+1=0\mod d$, and $b=-1\mod d$ cannot be limits 
in the sum \eqref{ademont} implying that
if a term of the form $j=-1 \mod d$ occurs
in the sum then so does $j+1=0 \mod d$.
Therefore, since these pairs of terms satisfy
$$
\omega^{j} \begin{bmatrix}
j-1 \\m-2
\end{bmatrix}_{\omega} + \omega^{j+1} \begin{bmatrix}
j \\m-2
\end{bmatrix}_{\omega} = \omega^{-1} 
\begin{pmatrix}
 n-1 \\m/d-1 
\end{pmatrix}+\frac{1-\omega^{-1}}{1-\omega}
\begin{pmatrix}
 n-1 \\m/d-1 
\end{pmatrix} =0 \, ,
$$
the proposition holds.
\end{proof}

\noindent
{\bf Acknowledgements.}
{\it We thank Dennis Stanton for his questions that initiated part of
this work and for several very helpful references.}


\begin{thebibliography}{33}
\bibitem{[GS]}
F. Garvan and D. Stanton, \emph{Sieved partition functions
and $q$-binomial coefficients}, Mathematics of computation.
{\bf 55} 191, 299--311 (1990).
\bibitem{[LLM]} L. Lapointe, A. Lascoux and J. Morse,
\emph{Tableau atoms and a
new Macdonald positivity conjecture}, 
Duke Math. J. {\bf 116}, 103--146 (2003).
\bibitem{[LM1]} 
L. Lapointe and J. Morse, \emph{Schur function analogs 
for a filtration
of the symmetric function space}, J. Comb. Th. A {\bf 101/2}, 191--224 (2003).
\bibitem{[LMrec]}
 L. Lapointe and J. Morse, \emph{Schur function identities,
their $t$-analogs, and $k$-Schur irreducibility}, Adv. in Math. {\bf 180},  
222--247 (2003). 
\bibitem{[LM3]} L. Lapointe and J. Morse, 
\emph{Tableaux on $k+1$-cores, reduced words for affine permutations, and 
$k$-Schur expansions}, http://www.math.miami.edu/\~{}morsej
\bibitem{[P]} R. Proctor, \emph{Solution of two difficult combinatorial
problems with linear algebra}, Amer. Math. Monthly {\bf 89}, 721--734 (1982).
\bibitem{[StEnu]} R.~P. Stanley, Enumerative Combinatorics, vol. 1, 
Cambridge
University Press, 1997.
\bibitem{[St]} R.~P. Stanley, {\it Weyl groups, the hard Lefschetz theorem,
and the Sperner property}, SIAM
J. Alg. Disc. Meth. {\bf 1}, 168--184 (1980).
\bibitem{[SW]} D. Stanton, and D. White, {\it Sieved $q$-binomial 
coefficients}, preprint.
\bibitem{[S]} J.~J. Sylvester, \emph{Proof of the hitherto undemonstrated 
fundamental theorem of invariants}, Collected Mathematical Papers,
vol. 3,  Chelsea, New York, 1973, 117--126.
\bibitem{[Ul]}  A. Ulyanov, \emph{Partitions with bounded part size and 
permutahedral tilings}, unpublished manuscript.
\bibitem{[WW]} S. Wangon, and H.S. Wilf, {\it When are subset sums
equidistributed modulo m?}, Electron. J. Combin. {\bf 1} (1994), R3. 
\bibitem{[DW]} D. Waugh, {\it Upper bounds in affine Weyl groups under
the weak order}, Order, {\bf 16}, 77-87 (1999).
\end{thebibliography}
\end{document}